\documentclass[final, 3p, times]{elsarticle}%
\usepackage{graphicx}
\usepackage{xcolor}
\usepackage{colortbl}
\usepackage{epstopdf}
\usepackage{amssymb}
\usepackage{amsthm}
\usepackage{amsmath}
\usepackage{color}
\usepackage{dsfont}
\usepackage{mathtools}
\usepackage{extarrows}
\usepackage{hyperref}
\usepackage{bm}
\usepackage{caption}
\usepackage{float}
\usepackage{verbatim}
\usepackage{epsfig}
\usepackage{multirow}
\usepackage[section]{placeins}
\usepackage{stmaryrd}
\usepackage{subfigure,tabularx,pdfpages}
\newcommand{\diff}{\,\mathrm{d}}

\newtheorem{lem}{Lemma}[section]
\newtheorem{thm}{Theorem}[section]

\newtheorem{exm}{Example}[section]
\numberwithin{equation}{section}
\biboptions{compress}

\begin{document}
\begin{frontmatter}
\title{Convergence analysis of exponential time differencing scheme for the nonlocal Cahn-Hilliard equation}

\tnotetext[label1]{The research of Dongling Wang is supported in part by the National Natural Science Foundation of China under grants 12271463 and Outstanding youth fund of department of education of Hunan province under grants 22B0173. The work of Danni Zhang is supported by Hunan Provincial Innovation Foundation For Postgraduate (No: CX20230625).
 \\ Declarations of interest: none.
}

\author[XTU]{Danni Zhang} 
\ead{zhdanni0807@smail.xtu.edu.cn;}
\author[XTU]{Dongling Wang\corref{mycorrespondingauthor}}
\ead{wdymath@xtu.edu.cn; ORCID 0000-0001-8509-2837}
\cortext[mycorrespondingauthor]{Corresponding author. }

\address[XTU]{School of Mathematics and Computational Science, Xiangtan University, Xiangtan, Hunan 411105, China}

\begin{abstract}
In this paper, we present a rigorous proof of the convergence of first order and second order exponential time differencing (ETD) schemes for solving the nonlocal Cahn-Hilliard (NCH) equation. The spatial discretization employs the Fourier spectral collocation method, while the time discretization is implemented using ETD-based multistep schemes. The absence of a higher-order diffusion term in the NCH equation poses a significant challenge to its convergence analysis. To tackle this, we introduce new error decomposition formulas and employ the higher-order consistency analysis.
These techniques enable us to establish the $\ell^\infty$ bound of numerical solutions under some natural constraints.
By treating the numerical solution as a perturbation of the exact solution,
we derive optimal convergence rates in $\ell^\infty(0,T;H_h^{-1})\cap \ell^2(0,T; \ell^2)$.
We conduct several numerical experiments to validate the accuracy and efficiency of the proposed schemes, including convergence tests and the observation of long-term coarsening dynamics.

\end{abstract}


\end{frontmatter}

\section{Introduction}\label{sec:Introd}

The Cahn-Hilliard (CH) equation was originally proposed in \cite{cahn} as one of the most typical phase field models which provides a macroscopic description of phase separation and microstructure evolution in binary alloy systems.
As a nonlocal variant of the classic CH equation, the Nonlocal Cahn-Hilliard (NCH) equation has attracted increasing attention and has been widely applied in diverse fields ranging from chemistry, material science and image processing \cite{3,4,5,7,Bates25}.
In this paper, we consider the following NCH equation with the periodic boundary conditions \cite{Bates25,Du1}:
\begin{equation}\label{equation:eq-(1.1)}
u_t=\Delta(\varepsilon^2\mathcal{L}u+f(u)),\quad (x,t)\in \Omega\times(0,T],
\end{equation}
where $\Omega=\prod_{i=1}^d(-X_i,X_i)$ is a rectangular domain in $\mathbb{R}^{d}(d=1,2,3)$ and $u=u(x,t)$ is an order parameter subject to the initial condition $u(x,0)=u_0(x)$, $T>0$ is the terminal time, $\varepsilon>0$ is an interfacial parameter.
The function $f(u)=F'(u)$ and $F(u)=\frac{1}{4}(u^2-1)^2$ is a double well function. The nonlocal linear operator $\mathcal{L}$ is defined by:
\begin{equation}\label{equation:eq-(1.2)}
\mathcal{L}: v(x) \longmapsto \int_{\Omega}J(x-y)(v(x)-v(y))\diff y,
\end{equation}
where $J$ is a nonnegative, $\Omega$-periodic radial kernel function and has a finite second moment in $\Omega$.

The linear operator $\mathcal{L}$ with the kernel function $J$ is self-adjoint and positive semi-definite.
Further, if $J$ is integrable, then $\mathcal{L}v=(J\ast 1)v-J\ast v$, where 
\begin{equation}
(J*v)(x)=\int_{\Omega}J(x-y)v(y)\diff y=\int_{\Omega}J(y)v(x-y)\diff y.
\end{equation}
We assume $\gamma_0:=\epsilon^2(J\ast 1)-1>0$ to ensure the NCH equation \eqref{equation:eq-(1.1)} is positive diffusive \cite{Bates1,Bates2}.

The NCH equation \eqref{equation:eq-(1.1)} can be view as the $H^{-1}$ gradient flow with respect to the free energy functional
\begin{equation}\label{equation:eq-(1.4)}
E(u)=\int_{\Omega}F(u(x))\diff x+\frac{\varepsilon^2}{2}(\mathcal{L}u,u),
\end{equation}
where $(\cdot,\cdot)$ denotes the standard $L^2$ inner product on $\Omega$.
For the smooth solutions of \eqref{equation:eq-(1.1)}, the total mass is conserved:
\begin{equation}\label{equation:eq-mass}
\frac{\diff}{\diff t} \int_\Omega u(x,t) \diff x \equiv 0.
\end{equation}
In this paper, we will only consider initial data $u_0$ with mean zero. Then, the fractional Laplacian operator $|\nabla|^s=(-\Delta)^{\frac{s}{2}}$ for $s<0$ is well defined \cite{li2016}.
Due to the gradient structure of \eqref{equation:eq-(1.1)}, the energy dissipation laws hold:
\begin{equation}
\frac{\diff}{\diff t}E(u)=-\||\nabla|^{-1}u_t\|_{L^2(\Omega)}^2 \leq 0.
\end{equation}

There have been many works on both theoretical analysis and numerical methods for the NCH equation.
In mathematical analysis, the well-posedness of the NCH equation equipped with a Neumann or Dirichlet boundary conditions were studied in \cite{Bates2,Bates1} with an integrable kernel function and in \cite{Guan1} claimed that the existence and uniqueness of the solution to the NCH equation subject to the periodic boundary condition may be established by using a similar technique.

For numerical methods, because of the energetic variational structure inherent in the nonlocal phase field model, an important fact is that the exact solution decrease the energy in time to NCH equation.
Therefore, it is essential to develop the so called energy stable numerical algorithms share this key property at the discrete time level.
Energy stable numerical schemes of the NCH equation including convex splitting schemes \cite{Guan1,Guan3}, exponential time differencing (ETD) schemes \cite{zhou2023,zhang2024,fu2024}, stabilized schemes \cite{Du1,li2023,li2024}, and so on.

The convergence of numerical solution to NCH equation are established in \cite{Guan1,Guan2, li1, li2023,li2024}.
In \cite{Guan1, Guan2}, convex splitting schemes were introduced for the NCH equation, demonstrating energy stability and convergence properties. These schemes handled the nonlinear term implicitly and the nonlocal term explicitly, leading to effective numerical solutions but nonlinear iterations were required.
Du et al. proposed energy-stable linear semi-implicit schemes in \cite{Du1}, utilizing stabilization techniques to avoid nonlinear iterations. The work in \cite{li1} presented the convergence in the discrete $H^{-1}$ norm and established the $\ell^\infty$ bound of the numerical solutions for the first-order stabilized linear semi-implicit scheme proposed in \cite{Du1}.
Furthermore, Li et al. \cite{li2023,li2024} developed two energy-stable and convergent second-order linear numerical schemes for the NCH equation. These schemes involved combining a modified Crank-Nicolson (CN) approximation with the Adams-Bashforth extrapolation and stability of the second-order backward differentiation formula (BDF2) for the time discretization, resulting in accurate and stable numerical solutions for the NCH equation.

However, the convergence analysis of ETD schemes for the NCH equation is more challenging and there are currently no results, one of the reasons being that
the lack of higher-order diffusion term and the Laplacian of nonlinear term.
The ETD schemes \cite{beylkin1998,cox2002,du2004,du2005,hochbruck2010}, which involve exact integration of the target equation followed by an explicit approximation of the temporal integral of the nonlinear term.
Exact evaluation of the linear terms makes the ETD schemes achieve high accuracy and satisfactory stability when dealing with stiff differential equations.
The first- and second-order ETD schemes were applied to the Nonlocal Allen-Cahn equation and have been analyzed energy stable and $\ell^\infty$ convergence \cite{Du2}, and further extended to a class of semi-linear parabolic equations \cite{Du3}.
Convergence and the $\ell^\infty$ bound of numerical solutions were proved of ETD schemes to classic CH equation \cite{li2019}.

Compare with the classic CH equation or Nonlocal Allen-Cahn equation, the NCH equation equipped with nonlocal diffusion operator and has more complicated  nonlinear term, so that the analytical techniques mentioned above are hardly applicable to the NCH equation.
In this paper, we establish the convergence of the fully discrete first-order ETD scheme (ETD1) and second-order ETD multistep scheme (ETD2) for the NCH equation which are developed our recent work \cite{zhang2024}.
Our analysis is mainly based on two important observations:
\begin{itemize}
\item
Unlike the standard error analysis based on the error equation \eqref{equation:eqe},
our analysis based on two new error decomposition formulas
\begin{align*}\label{eq:errd0}
&e^n:=\tilde{e}^n+\tau \mathcal{P}_N u_{\tau,1}+\tau^2 \mathcal{P}_N u_{\tau,2}\quad (\text{ETD1 scheme});\\ &e^n:=\tilde{e}^n+\tau^2 \mathcal{P}_N u_{\tau,3} \hspace{1.85cm} (\text{ETD2 scheme}).
\end{align*}
See the details in \eqref{eq:errd}.
This  formula allows us to establish the higher order consistency and obtain the $\ell^{\infty}$ bound of the numerical solutions.

\item
We adopt $(-\Delta_N)^{-1}\tilde{e}^{n+1}$ to test the error functions in \eqref{equation:eq-(3.19)} and \eqref{equation:eq-(3.44)}, instead of $\tilde{e}^{n+1}$ as in a previous work \cite{li2019}, where $(-\Delta_N)^{-1}$ is a spatial discrete operator defined in the Appendix A. This  method can be used to overcome the possible instability caused by the absence of high-order diffusion terms in the NCH equation.
\end{itemize}

Based on the above two points, we prove the convergence in the discrete $H^{-1}$ norm and obtain the $\ell^\infty$ bound of numerical solutions under some natural constraints on the space-time step sizes.

The rest of this paper is organized as follows.
The fully discrete ETD1 and ETD2 schemes are constructed in Section \ref{sec:ETD}, and some notations, definitions and useful lemmas are  also introduced.
In Section \ref{sec:Conver},  we prove the convergence by the higher-order consistency estimate. In addition, the $\ell^\infty$ bound of numerical solutions are also obtained.
In Section \ref{sec:Numer}, some numerical experiments are carried out to test the convergence in time level. The coarsening dynamics are simulated to show the long time behaviors of ETD2.
Some concluding remarks are given in Section \ref{sec:Conclu}.

\section{Fully discrete exponential time differencing schemes}\label{sec:ETD}
In this section, we present the fully discrete exponential time differencing schemes for the NCH equation, where the Fourier spectral collocation method is adopted for the spatial discretization.

We consider the square domain $\Omega=(-X,X)^2$ and define the spatial grid
$
\Omega_h=\{(x_i,y_j)=(-X+ih,-X+jh),\; 1\leq i,j\leq N\},
$
where the space step size $h=\frac{2X}{N}$ ($N$ is even).
We need two spaces. Let $\mathcal{M}_h$ be the space of all the $\Omega_h$-periodic grid functions.
That means if $f\in \mathcal{M}_h$, then $f$ is a discrete function defined on the discrete grid $\Omega_{h}$.
Let $\mathcal{M}_h^0$ be the space of the zero-mean grid functions. That means
if $f\in \mathcal{M}_h^0$, then $\overline{f}:=\frac{h^2}{4X^2}\sum_{(i, j)\in S_h} f_{i, j}=0$, where $S_h$ is the index set defined in the Appendix A.

 We denote $\Delta_N$ and $\mathcal{L}_N$ are Fourier spectral collocation space approximation operators of $\Delta$ and $\mathcal{L}$. These detailed definitions and some related properties are given in the Appendix A. For any $f, g\in \mathcal{M}_h$,
the discrete inner product $\langle\cdot,\cdot\rangle$, $\langle\cdot,\cdot\rangle_{-1,N}$ and norm $\|\cdot\|_2$, $\|\cdot\|_\infty$, $\|\cdot\|_{-1,N}$, and the discrete convolution $f\circledast \phi$ are defined in the Appendix A.

It is well known \cite{du2004} that a suitable linear operator splitting can improve the stability.
So we can introduce a stabilizing parameter $\kappa>0$ and define
\begin{equation}
L_h:=-\varepsilon^2\Delta_N \mathcal{L}_N-\kappa \Delta_N,\;\,\; f_\kappa(U):=f(U)-\kappa U.
\end{equation}
Now $L_h$ is self-adjoint and positive definite on $\mathcal{M}_h^0$. Then the space semi-discrete equation for \eqref{equation:eq-(1.1)} is to find $U:[0,T]\to \mathcal{M}_h^0$ such that
\begin{equation}\label{equation:eq-(2.9)}
\begin{split}
\left\{
\begin{split}
&\frac{\diff U}{\diff t}+L_hU=\Delta_N f_\kappa(U), \quad t\in(0,T], \\
&U(0)=U^0,
\end{split}
\right.
\end{split}
\end{equation}
where the initial value $U^0 \in \mathcal{M}_h^0$ is given.

Let $\{t_n=n\tau:0\leq n\leq N_t\}$, where $\tau=\frac{T}{N_t}$ is the time step size and $N_t$ is a given positive integer.
By using the property of the differentiation of matrix exponentials, the solution of the equation \eqref{equation:eq-(2.9)} satisfies
\begin{equation}\label{equation:eq-(2.100)}
U(t_{n+1})=e^{-\tau L_h}U(t_n)+\int_0^{\tau}e^{-(\tau-s)L_h}\Delta_N f_\kappa(U(t_n+s))\diff s.
\end{equation}
The key to construct the ETD schemes is to approximate the nonlinear term $f_\kappa(U(t_n+s))$ in \eqref{equation:eq-(2.100)} by polynomial interpolation, and then precisely integrate the interpolation.

The ETD1 scheme comes from approximating $f_\kappa(U(t_n+s))$ by the constant $f_\kappa(U(t_n))$, given by
\begin{equation}\label{equation:eq-(2.12)}
U^{n+1}=\phi_{-1}(\tau L_h)U^n+\tau \phi_{0}(\tau L_h)\Delta_Nf_\kappa(U^n),
\end{equation}
where $\phi_{-1}(a)=e^{-a},\phi_{0}(a)=\frac{1-e^{-a}}{a}, a\neq 0$. 
The ETD2 scheme is obtained by approximating $f_\kappa(U(t_n+s))$ by a linear extrapolation based on $f_\kappa(U(t_n))$ and $f_\kappa(U(t_{n-1}))$, given by
\begin{equation}\label{equation:eq-(2.13)}
U^{n+1}=\phi_{-1}(\tau L_h)U^n+\tau [(\phi_{0}+\phi_{1})(\tau L_h)\Delta_Nf_\kappa(U^n)-\phi_{1}(\tau L_h)\Delta_Nf_\kappa(U^{n-1})],
\end{equation}
where $\phi_{1}(a)=\frac{a-1+e^{-a}}{a^2}, a\neq 0.$

The ETD2 scheme \eqref{equation:eq-(2.13)} is a two-step algorithm, an accurate approximation for the value at $t_1$
is needed. Usually, $U^1$ can be computed by the ETD1.  But in this paper, we choose a higher-order approximation for $U^1$, to facilitate the higher-order consistency analysis presented in the later sections. For example, we can apply the second-order Runge-Kutta (RK2) method in first step, which yields a third-order temporal accuracy at $t_1$ for exact initial data $U^0$.

The following lemmas are useful in our following convergence analysis.
\begin{lem}[\cite{li2019}]\label{lem3}
\emph{(i)} For $a>0$, the following inequalities hold:
$
0<(1+a)\phi_{-1}(a)<1,\, 1<(1+a)\phi_0(a)<\frac{3}{2}, \, \frac{1}{2}<(1+a)\phi_1(a)<1$ and $0<(1+a)[\phi_0(a)-\phi_1(a)]<1.
$

\emph{(ii)}  If $0<s<\tau\leq 1$, then for any $a>0$, it holds
$0<(1+a \tau)e^{-a (\tau-s)}<1.$
\end{lem}
\begin{lem} [\cite{Du1}] \label{non}
The operator $\mathcal{L}_N$ has the following properties:

\emph{(i)}  the eigenvalues of $\mathcal{L}_N$ are $\lambda_{kl}=h^2(\hat{J}_{00}-\hat{J}_{kl})\geq 0,\, (k,l)\in \widehat{S}_h$;

\emph{(ii)} $\mathcal{L}_N$ commutes with $\Delta_N$ and is self-adjoint and positive semi-define;

\emph{(iii)}  for any $f\in\mathcal{M}_h$, we have $\mathcal{L}_Nf=(J\circledast 1)f-J\circledast f$.
\end{lem}
\begin{lem}[\cite{li1}]\label{lem1}
Suppose $J\in C_{per}^1(\Omega)$ and define its gird restriction via $J_{ij}:=J(x_i,y_j)$. Then for any $\phi,\psi\in \mathcal{M}_h$ and $\alpha>0$, we have
\begin{equation}
|\langle J\circledast\phi,\Delta_N \psi\rangle|\leq \alpha\|\phi\|_2^2+\frac{C_J}{\alpha}\|\nabla_N\psi\|_2^2,
\end{equation}
where $C_J$ is a positive constant that depends on $J$ and $\Omega$ but is independent of $h$.
\end{lem}

\section{Convergence analysis}\label{sec:Conver}
In this section, we analyze the convergence of the ETD1 scheme and ETD2 scheme, and give an optimal rate error estimates in $\ell^\infty(0,T;H_h^{-1})\cap \ell^2(0,T,\ell^2)$. To perform the convergence analysis, let's first introduce the concept and preliminary.

For a linear symmetric operator $Q: \mathcal{M}_h\to\mathcal{M}_h$, we define the norm $\interleave Q\interleave$ by the spectrum radius of $Q$, i.e., $\interleave Q\interleave =\max\{|\lambda|:\lambda\in\sigma(Q)\}$, where $\sigma(Q)$ is the set of all the eigenvalues of $Q$. It holds that
$
\|Qv\|_2\leq \interleave Q\interleave\cdot \|v\|_2,\;\forall v\in \mathcal{M}_h.
$

Denote by $u_e$ the exact solution of \eqref{equation:eq-(1.1)}. The existence and uniqueness of $u_e$ may be obtained by using the techniques adopted in \cite{Bates2,Bates1}, from which one can get the following estimate
\begin{equation}\label{equation:eq-(3.1)}
\|u_e\|_{L^\infty(0,T;L^\infty(\Omega))}+\|\nabla u_e\|_{L^\infty(0,T;L^\infty(\Omega))}\leq C, \quad \forall T>0.
\end{equation}
Moreover, if the initial data is sufficient regularity, we can assume that the exact solution has regularity as
\begin{equation}\label{eq:reg}
u_e\in\mathcal{R}:=H^4 \left(0,T;C^0_{per}(\bar{\Omega})\right)
\cap H^3\left( 0,T;C^2_{per}(\bar{\Omega}) \right)\cap L^{\infty} \left( 0,T;C^{m+2}_{per}(\bar{\Omega}) \right),\quad m\geq 3.
\end{equation}

Let $\mathcal{B}^K$ be the space of trigonometric polynomials of degree up to $K:=\frac{N}{2}$.
Let $\mathcal{P}_N: L^2(\Omega)\to \mathcal{B}^K$ be the $L^2$ orthogonal projection operator.
We define $u_N(\cdot, t):=\mathcal{P}_Nu_e(\cdot, t)$.
If $u_e\in L^\infty(0,T;H^\ell_{per}(\Omega))$ for some $\ell \in \mathbb{N}$, the following estimate is standard \cite{wang}
\begin{equation}\label{equation:eq-(3.2)}
\|u_N-u_e\|_{L^\infty(0,T;H^k(\Omega))}\leq Ch^{\ell-k}\|u_e\|_{L^\infty(0,T;H^\ell(\Omega))},\;\;\; \forall 0\leq k\leq \ell.
\end{equation}
By the orthogonality of Fourier projection $\mathcal{P}_N$, we have
$\int_{\Omega}u_N(\cdot,t_n)\diff x=\int_{\Omega}u_e(\cdot,t_n)\diff x$.
Noting that the exact solution $u_e$ is mass conservative, i.e., $\int_{\Omega}u_e(\cdot,t_n)\diff x=\int_{\Omega}u_e(\cdot,t_{n-1})\diff x$. Thus, we obtain
\begin{equation}\label{equation:eq-(3.3)}
\int_{\Omega}u_N(\cdot,t_n)\diff x=\int_{\Omega}u_N(\cdot,t_{n-1})\diff x,\quad \forall n\in\mathbb{N}.
\end{equation}
Since $u_N\in \mathcal{B}^K$, we have
$\int_{\Omega}u_N(\cdot,t_n)\diff x=h^2\sum_{(i,j)\in S_h} u_N(x_i,y_j,t_n)$, that means the rectangular quadrature rule holds exactly for all $u_N\in \mathcal{B}^K$.
The values of $u_N(t_n)$ at discrete grid points are still denoted as $u_N(t_n)$, i.e., $u_N(x_i, y_j, t_n)=\mathcal{P}_Nu_e(t_n)|_{i,j}$.
Recall~ $\overline{u_N}(t_n)=\frac{h^2}{4X^2}\sum_{(i, j)\in S_h} u_N(x_i,y_j,t_n),$
then from \eqref{equation:eq-(3.3)} we can get the mass conservative property of $u_N$ in the discrete average sense:
$\overline{u_N}(t_n)=\overline{u_N}(t_{n-1}).$

On the other hand, the solutions of the numerical schemes \eqref{equation:eq-(2.12)}-\eqref{equation:eq-(2.13)} are also mass conservative \cite{zhou2023}:
$
\overline{U^n}=\overline{U^{n-1}},\;\forall n\in\mathbb{N}.
$
We use  the mass conservative for the initial value $U^0=u_N(t_0)$,
and define the error grid function
\begin{equation}\label{equation:eqe}
e^n:=U^n-u_N(t_n),\quad \forall n\geq 0.
\end{equation}
We have $\overline{e^n}=0$, so $\|e^n\|_{-1,N}$ is well defined.
Under the regularity assumption \eqref{equation:eq-(3.1)}, for the projection functions $u_N(t_k):=\mathcal{P}_Nu_e(t_k)$, we have
\begin{equation}\label{eq:regu}
\max_{1\leq k\leq N_t}\|u_N(t_k)\|_{\infty}+\max_{1\leq k\leq N_t}\|\nabla_Nu_N(t_k)\|_{\infty}\leq C.
\end{equation}

The main result of this work is the following theorem.
\begin{thm}\label{thm:th1}

Assume that the solution of the NCH equation \eqref{equation:eq-(1.1)} satisfying the regularity class $\mathcal{R}$ presented in \eqref{eq:reg}. Also assume that $\tau$ and $h$ are small sufficiently and satisfy $\tau\leq Ch$ for some constant $C>0$. Let $B:=\max\limits_{1\leq k\leq N_t} \|u_N(t_k)\|_\infty+1$.

\emph{(i)} For ETD1 scheme, if the stabilizing parameter $\kappa\geq 2B^2-1$, we have
\begin{equation}\label{equation:eq1}
\|e^{n}\|_{-1,N}+ \left( \gamma_1\tau \sum_{k=1}^{n}\|e^k\|_2^2 \right)^{1/2}\leq C(\tau+h^m),\quad 0\leq n\leq N_t;
\end{equation}

\emph{(ii)} For ETD2 scheme, if the stabilizing parameter $\kappa\geq 3B^2-1$, we have
\begin{equation}\label{equation:eq2}
\|e^{n}\|_{-1,N}+ \left( \gamma_2\tau \sum_{k=1}^{n}\|e^k\|_2^2 \right)^{1/2}\leq C(\tau^2+h^m),\quad 0\leq n\leq N_t,
\end{equation}
where $C, \gamma_1, \gamma_2$ are some positive constant and are independent of $\tau$ and $h$.
\end{thm}

We provide a heuristic explanation of the main idea behind the proof.
The usual convergence analysis starts directly from the error equation $e^n:=U^n-u_N(t_n)$ defined in \eqref{equation:eqe}.
However, our here adopts a new strategy.
We construct two new error decomposition representation formulas
\begin{equation} \label{eq:errd}
\begin{split}
&e^n:=\tilde{e}^n+\tau \mathcal{P}_N u_{\tau,1}+\tau^2 \mathcal{P}_N u_{\tau,2}\quad (\text{ETD1 scheme}); \\
&e^n:=\tilde{e}^n+\tau^2 \mathcal{P}_N u_{\tau,3}\hspace{1.85cm} (\text{ETD2 scheme}),
\end{split}
\end{equation}
 where $\tilde{e}^n:=U^n-\tilde{u}(t_n)$ and $\tilde{u}(t_n)$ is the constructed approximate solution defined in \eqref{equation:eq-(3.6)} for the ETD1 scheme and in \eqref{equation:eq-(3.37)} for the ETD2 scheme, $u_{\tau,1}$, $u_{\tau,2}$ and $u_{\tau,3}$ are temporal correction functions.
 By doing so, we can obtain an $\mathcal{O}(\tau^3+h^m)$ truncation error.
 This higher consistency allows us to derive a higher order convergence estimate of the
modified error function $\tilde{e}^{n+1}$ as $\|\tilde{e}^{n+1}\|_{-1,N}= \|\tilde{u}(t_{n+1})-U^{n+1}\|_{-1,N}= \mathcal{O}(\tau^3+h^m)$ (see \eqref{equation:eq-(3.29)} and \eqref{equation:eq-(3.52)}), which in turn lead to $\|U^{n+1}\|_{\infty}<\infty$ under reasonable condition $\tau\le Ch$.
Due to the lack of higher-order diffusion term, we adopt $(-\Delta_N)^{-1}\tilde{e}^{n+1}$ to test the error function in \eqref{equation:eq-(3.19)} and \eqref{equation:eq-(3.44)}, instead of  $\tilde{e}^{n+1}$ as that in \cite{li2019}. Then, we estimate both sides of the error equation by using the result of the higher-order consistency analysis and induction argument, we can obtain an optimal convergence rate in $\ell^\infty(0,T;H_h^{-1})\cap \ell^2(0,T,\ell^2)$ for ${e}^{n+1}$ of the proposed two numerical schemes.

\subsection{Convergence analysis of ETD1}\label{ConETD1}

\subsubsection{Higher-order consistency analysis}

For a given $U^n$, the solution $U^{n+1}$ of the ETD1 scheme \eqref{equation:eq-(2.12)} can be given by $U^{n+1}=W(\tau)$ with the function $W: [0,\tau]\to\mathcal{M}_h^0$ determined by the following evolution equation:
\begin{equation}\label{equation:eq-(2.14)}
\begin{split}
\left\{
\begin{split}
&\frac{\diff W(s)}{\diff s}+L_hW(s)=\Delta_N f_\kappa(U^n), \quad s\in(0,\tau), \\
&W(0)=U^n.
\end{split}
\right.
\end{split}
\end{equation}
Then, for the continuous NCH equation \eqref{equation:eq-(1.1)}, we can give a similar expression as follows. For given $u_e(x,t_n)$, the solution $u_e(x,t_{n+1})$ can be obtained by $u_e(x,t_{n+1})=w(x,\tau)$ with the function $w(x,s)$ satisfying
\begin{equation}\label{equation:eq-(2.144)}
\begin{split}
\left\{
\begin{split}
&\frac{\partial w(s)}{\partial s}=\varepsilon^2 \Delta \mathcal{L}w+\kappa \Delta w+\Delta f_\kappa(w), \quad x\in \Omega,  s\in(0,\tau), \\
&w(x,0)=u_e(x,t_n).
\end{split}
\right.
\end{split}
\end{equation}

Denote the function $W_N(s)=\mathcal{P}_Nw(s)$. Then $W_N(s)$ satisfies the following equation:
\begin{equation}\label{equation:eq-(3.88)}
\frac{\partial W_N(s)}{\partial s}=\varepsilon^2 \Delta \mathcal{L}W_N(s)+\kappa\Delta W_N(s)+\Delta\mathcal{P}_N(u_e(t_n+s)^3)-\Delta W_N(s)-\kappa\Delta W_N(s),
\end{equation}
and $W_N(s)$ also solves the discrete equation \eqref{equation:eq-(2.14)} with the local truncation error $R_{h\tau}^{(1)}(s)$
\begin{equation}\label{equation:eq-(3.888)}
\begin{split}
\left\{
\begin{split}
&\frac{\diff W_N(s)}{\diff s}+L_hW_N(s)=\Delta_N f_\kappa(u_N(t_n))+R_{h\tau}^{(1)}(s), \quad s\in(0,\tau), \\
&W_N(0)=u_N(t_n).
\end{split}
\right.
\end{split}
\end{equation}
By comparing equations  \eqref{equation:eq-(3.88)} and \eqref{equation:eq-(3.888)}, we can get
$R_{h\tau}^{(1)}(s)=R_h^{(1)}+R_\tau^{(1)},$
where
\begin{align*}
&R_h^{(1)}=(\varepsilon^2 \Delta \mathcal{L}-\varepsilon^2 \Delta_N \mathcal{L}_N)u_N(t_n+s)
+(\kappa \Delta -\kappa \Delta_N) u_N(t_n+s)\\
&\hspace{1cm}+\Delta(\mathcal{P}_N(u_e(t_n+s)^3)-u_N(t_n+s)^3)
+\Delta f_\kappa(u_N(t_n+s))-\Delta_N f_\kappa(u_N(t_n+s)),\\
&R_\tau^{(1)}=\Delta_N f_\kappa(u_N(t_n+s))-\Delta_N f_\kappa(u_N(t_n)).
\end{align*}
By the standard Fourier spectral approximation, $R_h^{(1)}$ has the bound $Ch^m$. By Taylor expansion with the integral remainder, $R_\tau^{(1)}$ can be bounded by $C\tau$. Thus, we have
\[
\sup\limits_{s\in(0,\tau)}\|R_{h\tau}^{(1)}(s)\|_{-1,N}\leq C(\tau+h^m).
\]

Since the local truncation error only has first order accuracy in time, it is not enough to derive the boundedness of numerical solution, i.e., $\|U^{n+1}\|_{\infty}<\infty$.
To remedy this, we introduce the auxiliary profile
\begin{equation}\label{equation:eq-(3.6)}
\tilde{u}=u_N+\tau \mathcal{P}_N u_{\tau,1}+\tau^2\mathcal{P}_N u_{\tau,2},
\end{equation}
where $u_{\tau,1}$ and $u_{\tau,2}$ are two temporal correction functions will be constructed later.
 By doing so, we can obtain a higher $\mathcal{O}(\tau^3+h^m)$ consistency.

Let's construct the correction function $u_{\tau,1}$ first.
By Fourier projection from continuous NCH equation \eqref{equation:eq-(1.1)} and carry out the space discretization, we have $u_N: (t_n,t_{n+1}] \to\mathcal{M}_h^0$ satisfying
\[
\frac{\diff u_N}{\diff t}+L_hu_N=\Delta_N f_\kappa(u_N(t))+\mathcal{O}(h^m),\quad  t\in(t_n,t_{n+1}].
\]
Expansing the nonlinear term $f_\kappa(u_N(t))$ for $t\in(t_n,t_{n+1}]$ at
$t_n$ to obtain $f_\kappa(u_N(t))=f_\kappa(u_N(t_{n}))+\tau g_{1}+\mathcal{O}(\tau^2)$, where $g_{1}$ depends only on $f_\kappa$ and $u_N$.
We get that
\begin{equation}\label{equation:eq-(3.8)}
\frac{\diff u_N}{\diff t}+L_hu_N=\Delta_N f_\kappa(u_N(t_n))+\tau g_1+\mathcal{O}(\tau^2)+\mathcal{O}(h^m),\quad  t\in(t_n,t_{n+1}].
\end{equation}
By mass conservative and the periodic boundary condition, we also have $\overline{g_1}=0$.

The first order temporal correction function $u_{\tau,1}$ is given by solving the following ODE system:
\begin{equation}\label{equation:eq-(3.9)}
\begin{split}
\left\{
\begin{split}
&\frac{\diff u_{\tau,1}}{\diff t}+L_hu_{\tau,1}=\Delta_N \left[ f'_\kappa(u_N(t_n)) u_{\tau,1}(t_n) \right]-g_1,\quad t\in(t_n,t_{n+1}],\\
&u_{\tau,1}(0)\equiv 0,
\end{split}
\right.
\end{split}
\end{equation}
The existence and uniqueness of the solution of this linear system is standard and it depends only on $u_N$.

Let $\tilde{u}_1=u_N+\tau \mathcal{P}_Nu_{\tau,1}$.
Then it follows from the equations \eqref{equation:eq-(3.8)} and \eqref{equation:eq-(3.9)} that
\begin{equation}\label{equation:eq-(3.10)}
\frac{\diff \tilde{u}_1}{\diff t}+L_h\tilde{u}_1=\Delta_N f_\kappa(\tilde{u}_1(t_n))+\tau^2g_2+\mathcal{O}(\tau^3+h^m),\quad t\in(t_n,t_{n+1}],
\end{equation}
where we have used the following expansion
\begin{align*}
f_\kappa(\tilde{u}_1(t_n))&=f_\kappa(u_N(t_n))+\tau f_\kappa'(u_N(t_n))\mathcal{P}_Nu_{\tau,1}(t_n)+\mathcal{O}(\tau^2)\\
&=f_\kappa(u_N(t_n))+\tau \mathcal{P}_N(f_\kappa'(u_N(t_n))\mathcal{P}_Nu_{\tau,1}(t_n))+\mathcal{O}(\tau^2+h^m),
\end{align*}
and $g_2$ defined in \eqref{equation:eq-(3.10)} depends only on $f_\kappa$ and $u_N$. We also have $\overline{g_2} =0$.
Similarly, we can get the second order temporal correction function $u_{\tau,2}$, which is the solution of the following linear ODE system:
\begin{equation}\label{equation:eq-(3.11)}
\begin{split}
\left\{
\begin{split}
&\frac{\diff u_{\tau,2}}{\diff t}+L_hu_{\tau,2}=\Delta_N[f_\kappa'(u_N(t_n))u_{\tau,2}(t_n)]-g_{2},\quad t\in(t_n,t_{n+1}], \\
&u_{\tau,2}(0)\equiv 0,
\end{split}
\right.
\end{split}
\end{equation}
The unique solution only dependent on $u_N$. Define a grid function $\tilde{u}=\tilde{u}_1+\tau^2 \mathcal{P}_Nu_{\tau,2}$. Combination of \eqref{equation:eq-(3.10)} and a Fourier projection of \eqref{equation:eq-(3.11)}  leads to
\begin{equation}\label{equation:eq-(3.12)}
\frac{\diff \tilde{u}}{\diff t}+L_h\tilde{u}=\Delta_N f_\kappa(\tilde{u}(t_n))+\mathcal{O}(\tau^3+h^m),\quad t\in(t_n,t_{n+1}],
\end{equation}
where we have used the fact
\begin{align*}
f_\kappa(\tilde{u}(t_n))&=f_\kappa(\tilde{u}_1(t_n))+\tau^2 f_\kappa'(u_N(t_n))\mathcal{P}_Nu_{\tau,2}(t_n)+\mathcal{O}(\tau^3)\\
&=f_\kappa(\tilde{u}_1(t_n))+\tau^2 \mathcal{P}_N(f_\kappa'(u_N(t_n))\mathcal{P}_Nu_{\tau,2}(t_n))+\mathcal{O}(\tau^3+h^m).
\end{align*}
We see that $\overline{\tilde{u}}=0$.
Thus, we have completed the construction of $\tilde{u}$.

Note that $\|\tilde{u}\|_\infty$ is also bounded. In fact, it follows from $\|u_{\tau,1}\|_\infty\leq C$ and $\|u_{\tau,2}\|_\infty\leq C$ that
$\|\tilde{u}-u_N\|_\infty\leq C\tau.$
In particular, when $\tau$ is sufficiently small we have
$\|\tilde{u}-u_N\|_\infty\leq C\tau\leq \frac{1}{2}$,
which gives
$\|\tilde{u}\|_\infty\leq\|u_N\|_\infty+ \|\tilde{u}-u_N\|_\infty\leq \|u_N\|_\infty+\frac{1}{2}\leq B.$

\subsubsection{Convergence analysis in $\ell^\infty \left( 0,T;H_h^{-1} \right) \cap \ell^2 \left( 0,T; \ell^2 \right)$}

Let $\tilde{W}(s)=\tilde{u}(t_n+s)$ for $s\in[0,\tau]$. Then, by consistency, $\tilde{W}(s)$ satisfies
\begin{equation}\label{equation:eq-(3.17)}
\begin{split}
\left\{
\begin{split}
&\frac{\diff \tilde{W}(s)}{\diff s}+L_h\tilde{W}(s)=\Delta_N f_\kappa(\tilde{u}(t_n))+\tilde{R}_{h\tau}^{(1)}(s), \quad  s\in(0,\tau), \\
&\tilde{W}(0)=\tilde{u}(t_n),
\end{split}
\right.
\end{split}
\end{equation}
where the truncation error $\tilde{R}_{h\tau}^{(1)}(s)$ satisfies
\begin{equation}\label{equation:eq-(3.181)}
\sup\limits_{s\in(0,\tau)}\|\tilde{R}_{h\tau}^{(1)}(s)\|_{-1,N}\leq C(\tau^3+h^m),
\end{equation}
Subtracting \eqref{equation:eq-(3.17)} from \eqref{equation:eq-(2.14)} and putting $\tilde{e}(s):=W(s)-\tilde{W}(s)$ yield
\begin{equation}\label{equation:eq-(3.18)}
\begin{split}
\left\{
\begin{split}
&\frac{\diff \tilde{e}(s)}{\diff s}+L_h\tilde{e}(s)=\Delta_N f_\kappa(U^n)-\Delta_N f_\kappa(\tilde{u}(t_n))-\tilde{R}_{h\tau}^{(1)}(s), \quad s\in(0,\tau), \\
&\tilde{e}(0)=U^n-\tilde{u}(t_n)=:\tilde{e}^n\in \mathcal{M}_h^0.
\end{split}
\right.
\end{split}
\end{equation}
Hence, $\tilde{e}^{n+1}:=\tilde{e}(\tau)$ satisfies that
\begin{equation}\label{equation:eq-(3.19)}
\tilde{e}^{n+1}=\phi_{-1}(\tau L_h)\tilde{e}^n+\tau\phi_{0}(\tau L_h)[\Delta_N f_\kappa(U^n)-\Delta_N f_\kappa(\tilde{u}(t_n))]-\int_0^{\tau}\tilde{e}^{-(\tau-s)L_h}\tilde{R}_{h\tau}^{(1)}(s) \diff s.
\end{equation}
Acting $I+\tau L_h$ on \eqref{equation:eq-(3.19)} and taking the discrete $\ell^2$ inner product with $(-\Delta_N)^{-1}\tilde{e}^{n+1}$ yield
\begin{equation}\label{equation:eq-(3.20)}
\|\tilde{e}^{n+1}\|_{-1,N}^2+\tau \varepsilon^2\langle \mathcal{L}_N \tilde{e}^{n+1},\tilde{e}^{n+1}\rangle+\tau \kappa\|\tilde{e}^{n+1}\|_2^2=R_1,
\end{equation}
where
\begin{align}\label{eq:R}
R_1=&\left\langle(I+\tau L_h)\phi_{-1}(\tau L_h)\tilde{e}^n,(-\Delta_N)^{-1}\tilde{e}^{n+1} \right\rangle\nonumber\\
&+\tau \left\langle(I+\tau L_h)\phi_{0}(\tau L_h)\Delta_N[f_\kappa(U^n)-f_\kappa(\tilde{u}(t_n))],(-\Delta_N)^{-1}\tilde{e}^{n+1} \right\rangle\nonumber\\
&-\int_0^{\tau}\left\langle(I+\tau L_h)\tilde{e}^{-(\tau-s)L_h}\tilde{R}_{h\tau}^{(1)}(s),(-\Delta_N)^{-1}\tilde{e}^{n+1}\right\rangle\diff s.
\end{align}

Next, we use induction argument to prove the convergence.
We note that numerical error function $\|\tilde{e}^0\|_\infty=0<\frac{1}{2}$ at $t=0$. Suppose $\|\tilde{e}^n\|_\infty\leq \frac{1}{2}$ at $t_n$. Thus, the $\ell^\infty$ bound for the numerical solutions at $t_n$ becomes available:
\[
\|U^n\|_\infty=\|\tilde{u}(t_n)+\tilde{e}^n\|_\infty \leq \|\tilde{u}(t_n)\|_\infty +\|\tilde{e}^n\|_\infty\leq \|u_N(t_n)\|_\infty+\frac{1}{2}+\frac{1}{2}\leq B.
\]

We now estimate the three terms on the right-hand side given in \eqref{eq:R}.
For the first linear term, a direct application of Cauchy inequality and Lemma \ref{lem3} give
\begin{align}\label{equation:eq-(3.23)}
\left\langle(I+\tau L_h)\phi_{-1}(\tau L_h)\tilde{e}^n,(-\Delta_N)^{-1}\tilde{e}^{n+1}\right\rangle&=\left\langle(I+\tau L_h)\phi_{-1}(\tau L_h)(-\Delta_N)^{-\frac{1}{2}}\tilde{e}^n,(-\Delta_N)^{-\frac{1}{2}}\tilde{e}^{n+1}\right\rangle\nonumber \\
&\leq \|(I+\tau L_h)\phi_{-1}(\tau L_h)(-\Delta_N)^{-\frac{1}{2}}\tilde{e}^n\|_2\|(-\Delta_N)^{-\frac{1}{2}}\tilde{e}^{n+1}\|_2\nonumber \\
&\leq \frac{1}{2}\|\tilde{e}^n\|_{-1,N}^2+\frac{1}{2}\|\tilde{e}^{n+1}\|_{-1,N}^2.
\end{align}

By noting that $\|U^n\|_\infty\leq B, \|\tilde{u}\|_\infty\leq B$, we have
\[
\|f_\kappa(U^n)-f_\kappa(\tilde{u}(t_n))\|_2=\|f_\kappa'(\xi_n)(U^n-\tilde{u}(t_n))\|_2\leq \|f_\kappa'(\xi_n)\|_\infty\|\tilde{e}^n\|_2\leq |3B^2-\kappa-1|\|\tilde{e}^n\|_2:= K\|\tilde{e}^n\|_2,
\]
where $\xi_n$ between $U^n$ and $\tilde{u}(t_n)$ and $K=|3B^2-\kappa-1|$.
Then, for the second term in \eqref{eq:R}, by using Lemma \ref{lem3} and the above estimate, we have
\begin{align}\label{equation:eq-(3.21)}
&\tau \left\langle(I+\tau L_h)\phi_{0}(\tau L_h)\Delta_N[f_\kappa(U^n)-f_\kappa(\tilde{u}(t_n))],(-\Delta_N)^{-1}\tilde{e}^{n+1} \right\rangle\nonumber\\
 &\leq \tau\interleave (I+\tau L_h)\phi_{0}(\tau L_h)\interleave\|f_\kappa(U^n)-f_\kappa(\tilde{u}(t_n))\|_2\|\tilde{e}^{n+1}\|_2\nonumber\\
 &< 2\tau K\|\tilde{e}^n\|_2 \|\tilde{e}^{n+1}\|_2\leq \tau K\|\tilde{e}^n\|_2^2+\tau K\|\tilde{e}^{n+1}\|_2^2.
\end{align}

For the last term in \eqref{eq:R}, we have the following estimate:
\begin{align}\label{equation:eq-(3.22)}
&-\int_0^{\tau}\left\langle(I+\tau L_h)\tilde{e}^{-(\tau-s)L_h}\tilde{R}_{h\tau}^{(1)}(s),(-\Delta_N)^{-1}\tilde{e}^{n+1}\right\rangle\diff s \nonumber\\
&\leq \int_0^{\tau}\interleave(I+\tau L_h)\tilde{e}^{-(\tau-s)L_h}\interleave \diff s \sup_{t\in(0,\tau)}\|\tilde{R}_{h\tau}^{(1)}(t)\|_{-1,N}\|\tilde{e}^{n+1}\|_{-1,N}\nonumber\\
&\leq \frac{\tau}{2}\|\tilde{e}^{n+1}\|_{-1,N}^2+\frac{\tau}{2}\sup_{t\in(0,\tau)}\|\tilde{R}_{h\tau}^{(1)}(t)\|_{-1,N}^2.
\end{align}

For the nonlocal linear term given in \eqref{equation:eq-(3.20)}, by appling Lemma \ref{non} and Lemma \ref{lem1}, we obtain
\begin{align}\label{equation:eq-(3.24)}
-\tau \varepsilon^2\langle \mathcal{L}_N \tilde{e}^{n+1},\tilde{e}^{n+1}\rangle&=-\tau \varepsilon^2\left\langle (J\circledast1) \tilde{e}^{n+1}-J\circledast \tilde{e}^{n+1},\tilde{e}^{n+1}\right\rangle\nonumber\\
&=-\tau \varepsilon^2(J\circledast1)\|\tilde{e}^{n+1}\|_2^2+\tau \varepsilon^2\langle J\circledast \tilde{e}^{n+1},\tilde{e}^{n+1} \rangle \rangle\nonumber\\
&=-\tau \varepsilon^2(J\circledast1)\|\tilde{e}^{n+1}\|_2^2-\tau \varepsilon^2\left\langle J\circledast \tilde{e}^{n+1},\Delta_N ((-\Delta_N)^{-1}\tilde{e}^{n+1}) \right\rangle\nonumber\\
&\leq-\tau \varepsilon^2(J\circledast1)\|\tilde{e}^{n+1}\|_2^2+\tau \gamma_0 \|\tilde{e}^{n+1}\|_2^2 +\tau \frac{C_1}{\gamma_0}\|\tilde{e}^{n+1}\|_{-1,N}^2,
\end{align}
where $C_1$ depends on $C_J$ and $\varepsilon$.

Therefore, a substitution of \eqref{equation:eq-(3.21)}-\eqref{equation:eq-(3.22)} and \eqref{equation:eq-(3.24)} into \eqref{equation:eq-(3.20)}, and recall the define of $\gamma_0$, we get
\begin{align}
&\frac{1}{2}(\|\tilde{e}^{n+1}\|_{-1,N}^2-\|\tilde{e}^n\|_{-1,N}^2)+(\kappa+1)\tau\|\tilde{e}^{n+1}\|_2^2\nonumber\\
&\leq \tau K \|\tilde{e}^n\|_2^2+\tau K\|\tilde{e}^{n+1}\|_2^2+(\frac{1}{2}+\frac{C_1}{\gamma_0})\tau \|\tilde{e}^{n+1}\|_{-1,N}^2+\frac{\tau}{2}\sup_{t\in(0,\tau)}\|\tilde{R}_{h\tau}^{(1)}(t)\|_{-1,N}^2.
\end{align}
Summing the above inequality from 0 to $n$ leads to
\begin{equation*}
\|\tilde{e}^{n+1}\|_{-1,N}^2+2(\kappa+1)\tau \sum_{k=0}^{n+1}\|\tilde{e}^k\|_2^2
\leq (1+\frac{2C_1}{\gamma_0})\tau \sum_{k=0}^{n+1}\|\tilde{e}^k\|_{-1,N}^2
+4\tau K\sum_{k=0}^{n+1}\|\tilde{e}^k\|_2^2+\tau\sum_{k=0}^n\sup_{t\in(0,\tau)}\|\tilde{R}_{h\tau}^{(1)}(t)\|_{-1,N}^2.
\end{equation*}
Subsequently, assuming that $\tau<(1+\frac{2C_1}{\gamma_0})^{-1}$ and putting $\gamma_1:=2(\kappa+1)-4K\geq 0$, we can apply  the discrete Gronwall's inequality to obtain the following convergence estimate:
\begin{equation}\label{equation:eq-(3.29)}
\|\tilde{e}^{n+1}\|_{-1,N}+ \left( \gamma_1\tau \sum_{k=1}^{n+1}\|\tilde{e}^k\|_2^2 \right)^{1/2}\leq C^*(\tau^3+h^m),
\end{equation}
where $C^*$ depends on $C_1, \gamma_0, T$.
By using the inverse inequality to \eqref{equation:eq-(3.29)}, the following estimate holds
\begin{align}\label{equation:eq-(3.291)}
\|\tilde{e}^{n+1}\|_\infty&\leq \frac{C\|\tilde{e}^{n+1}\|_{-1,N}}{h^2}\leq\frac{CC^*(\tau^3+h^m)}{h^2}\leq \frac{C'C^*(h^3+h^m)}{h^2}\leq\frac{C_2C^*h^3}{h^2}=C_2C^*h\leq \frac{1}{2},
\end{align}
provided that $h\leq \frac{1}{2C_2C^*}$, where we assumed $\tau\leq Ch$ and the fact that $m\geq 3$. Then, we have 
\begin{equation}\label{equation:eq-(3.292)}
\|U^{n+1}\|_\infty\leq \|\tilde{u}(t_{n+1})+\tilde{e}^{n+1}\|_\infty\leq B.
\end{equation}
This completes the error estimate for error function $\tilde{e}^{n+1}$.

Now, by the definition of the function $\tilde{u}=u_N+\tau \mathcal{P}_N u_{\tau,1}+ \tau^2 \mathcal{P}_N u_{\tau,2}$, we have
\begin{align}
\|e^{n+1}\|_{-1,N}&=\|U^{n+1}-u_N(t_{n+1})\|_{-1,N}\nonumber\\
&\leq \|\tilde{e}^{n+1}\|_{-1,N}+\|\tilde{u}(t_{n+1})-u_N(t_{n+1})\|_{-1,N}\nonumber\\
&\leq \|\tilde{e}^{n+1}\|_{-1,N}+\|\tau \mathcal{P}_N (u_{\tau,1}(t_{n+1}))+ \tau^2 \mathcal{P}_N (u_{\tau,2}(t_{n+1}))\|_{-1,N}\nonumber\\
&\leq \|\tilde{e}^{n+1}\|_{-1,N}+\tau\|\mathcal{P}_N (u_{\tau,1}(t_{n+1}))\|_{-1,N}+\tau^2\|\mathcal{P}_N (u_{\tau,2}(t_{n+1}))\|_{-1,N}.
\end{align}
In a similar manner,
\begin{equation}
\|e^{n+1}\|_2\leq \|\tilde{e}^{n+1}\|_2+\tau\|\mathcal{P}_N (u_{\tau,1}(t_{n+1}))\|_2+\tau^2\|\mathcal{P}_N (u_{\tau,2}(t_{n+1}))\|_2.
\end{equation}
Then, thanks to the above observations, we can conclude the error estimate \eqref{equation:eq1} from \eqref{equation:eq-(3.29)} and the fact that the boundedness of $u_{\tau,1}, u_{\tau,2}$. This competes the proof of Theorem \ref{thm:th1} (i).

\subsection{Convergence analysis of  ETD2}
Similar to the proof of convergence of ETD1 scheme in subsection \ref{ConETD1},
we next prove the convergence of the ETD2 scheme.
\subsubsection{Higher-order consistency analysis}

Similar to the higher-order consistency analysis of the ETD1 scheme, the Fourier projection solution $u_N$ satisfies the following equation:
\begin{equation}\label{equation:eq-(3.35)}
\frac{\diff u_N}{\diff t}+L_hu_N= \left( 1+\frac{t-t_n}{\tau} \right) \Delta_N f_\kappa(u_N(t_n))-\frac{t-t_n}{\tau}\Delta_N f_\kappa (u_N(t_{n-1}))+\tau^2 g_3+\mathcal{O}(\tau^3+h^m),\quad  t\in(t_n,t_{n+1}],
\end{equation}
where the function $g_3$ depends only on $f_\kappa$ and $u_N$ and we have $\overline{g_3} =0$.
Define the auxiliary profile
\begin{equation}\label{equation:eq-(3.37)}
\tilde{u}=u_N+\tau^2\mathcal{P}_Nu_{\tau,3},
\end{equation}
where the temporal correction function $u_{\tau,3}$ solves the following linear ODE system:
\begin{align}\label{equation:eq-(3.36)}
\frac{\diff u_{\tau,3}}{\diff t}+L_hu_{\tau,3}=&\left( 1+\frac{t-t_n}{\tau} \right) \Delta_N[f_\kappa'(u_N(t_n))u_{\tau,3}(t_n)]\nonumber \\
&-\frac{t-t_n}{\tau}\Delta_N [f_\kappa' (u_N(t_{n-1}))u_{\tau,3}(t_{n-1})]-g_3,\quad t\in(t_n,t_{n+1}],
\end{align}
subject to the zero initial value.
Note that the solution of \eqref{equation:eq-(3.36)} exists and is unique and the solution is smooth enough.
Combination of \eqref{equation:eq-(3.35)} and \eqref{equation:eq-(3.36)} gives the higher-order consistency estimate:
\begin{equation}\label{equation:eq-(3.38)}
\frac{\diff \tilde{u}}{\diff t}+L_h\tilde{u}=\left( 1+\frac{t-t_n}{\tau} \right) \Delta_N f_\kappa(\tilde{u}(t_n))-\frac{t-t_n}{\tau}\Delta_N f_\kappa (\tilde{u}(t_{n-1}))+\mathcal{O}(\tau^3+h^m),\quad t\in(t_n,t_{n+1}],
\end{equation}
where we have used the estimate
\begin{align*}
f_\kappa(\tilde{u}(t_k))&=f_\kappa(u_N(t_k))+\tau^2 f_\kappa'(u_N(t_k))\mathcal{P}_Nu_{\tau,3}(t_k)+\mathcal{O}(\tau^4)\\
&=f_\kappa(u_N(t_k))+\tau^2 \mathcal{P}_N(f_\kappa'(u_N(t_k))\mathcal{P}_Nu_{\tau,3}(t_k))+\mathcal{O}(\tau^4+h^m),\quad k=n,n-1.
\end{align*}

We also have $\overline{\tilde{u}}=0$ and the boundedness of $\|\tilde{u}\|_\infty$. In fact, by noting $\|u_{\tau,3}\|_\infty\leq C$ and if $\tau$ is small sufficiently, we have
$
\|\tilde{u}\|_\infty\leq\|u_N\|_\infty+ \|\tilde{u}-u_N\|_\infty\leq \|u_N\|_\infty+\frac{1}{2}\leq B.
$
And if we adopt the RK2 numerical algorithm for the estimate of $U^1$, it holds that
$
U^1-\tilde{u}(t_1)=\mathcal{O}(\tau^3+h^m).
$

\subsubsection{Convergence analysis in $\ell^\infty(0,T;H_h^{-1})\cap \ell^2(0,T,\ell^2)$}
For given $U^n, U^{n-1}\in \mathcal{M}_h^0$, the solution $U^{n+1}$ of the ETD2 scheme \eqref{equation:eq-(2.13)} is actually given by $U^{n+1}=W(\tau)$ with the function $W:[0,\tau]\to\mathcal{M}_h^0$ determined by the evolution equation:

\begin{equation}\label{equation:eq-(2.15)}
\begin{split}
\left\{
\begin{split}
&\frac{\diff W(s)}{\diff s}+L_hW(s)=\left( 1+\frac{s}{\tau} \right)\Delta_N f_\kappa(U^n)-\frac{s}{\tau}\Delta_N f_\kappa (U^{n-1}), \quad s\in(0,\tau), \\
&W(0)=U^n.
\end{split}
\right.
\end{split}
\end{equation}

The function $\tilde{W}(s)=\tilde{u}(t_n+s),\; s\in[0,\tau]$ satisfies \eqref{equation:eq-(2.15)} with the truncated error $\tilde{R}_{h\tau}^{(2)}(s)$
\begin{equation}\label{equation:eq-(3.42)}
\begin{split}
\left\{
\begin{split}
&\frac{\diff \tilde{W}(s)}{\diff s}+L_h\tilde{W}(s)=\left( 1+\frac{s}{\tau} \right) \Delta_N f_\kappa(\tilde{u}(t_n))-\frac{s}{\tau}\Delta_N f_\kappa(\tilde{u}(t_{n-1}))+\tilde{R}_{h\tau}^{(2)}(s), \quad s\in(0,\tau), \\
&\tilde{W}(0)=\tilde{u}(t_n),
\end{split}
\right.
\end{split}
\end{equation}
where $\tilde{R}_{h\tau}^{(2)}(s)$ satisfies
\[
\sup\limits_{s\in(0,\tau)}\|\tilde{R}_{h\tau}^{(2)}(s)\|_{-1,N}\leq C(\tau^3+h^m).
\]
Let $\tilde{e}(s):=W(s)-\tilde{W}(s)$. The difference between \eqref{equation:eq-(2.15)} and \eqref{equation:eq-(3.42)} gives
\begin{equation}\label{equation:eq-(3.43)}
\begin{split}
\left\{
\begin{split}
&\frac{\diff \tilde{e}(s)}{\diff s}+L_h\tilde{e}(s)=\left( 1+\frac{s}{\tau} \right) [\Delta_N f_\kappa(U^n)-\Delta_N f_\kappa(\tilde{u}(t_n))]\\
&\hspace{2.6cm}-\frac{s}{\tau}[\Delta_N f_\kappa(U^{n-1})-\Delta_N f_\kappa(\tilde{u}(t_{n-1}))]-\tilde{R}_{h\tau}^{(2)}(s), \quad s\in(0,\tau), \\
&\tilde{e}(0)=U^n-\tilde{u}(t_n):=\tilde{e}^n\in\mathcal{M}_h^0.
\end{split}
\right.
\end{split}
\end{equation}
Thus, the solution $\tilde{e}^{n+1}:=\tilde{e}(\tau)$ satisfies that
\begin{align}\label{equation:eq-(3.44)}
\tilde{e}^{n+1}=&\phi_{-1}(\tau L_h)\tilde{e}^n+\tau(\phi_0+\phi_1)(\tau L_h)\left[\Delta_N f_\kappa(U^n)-\Delta_N f_\kappa(\tilde{u}(t_n))\right]\nonumber\\
&-\tau \phi_1(\tau L_h)\left[\Delta_N f_\kappa(U^{n-1})-\Delta_N f_\kappa(\tilde{u}(t_{n-1}))\right]-\int_0^{\tau}\tilde{e}^{-(\tau-s)L_h}\tilde{R}_{h\tau}^{(2)}(s)\diff s.
\end{align}
Acting $I+\tau L_h$ on both sides of \eqref{equation:eq-(3.44)} and taking $\ell^2$ inner product of the resulted equation with $(-\Delta_N)^{-1}\tilde{e}^{n+1}$ yield
\begin{equation}\label{equation:eq-(3.45)}
\|\tilde{e}^{n+1}\|_{-1,N}^2+\tau \varepsilon^2\langle \mathcal{L}_N \tilde{e}^{n+1},\tilde{e}^{n+1}\rangle+\tau \kappa\|\tilde{e}^{n+1}\|_2^2=R_2,
\end{equation}
where
\begin{align*}
R_2=&\left\langle(I+\tau L_h)\phi_{-1}(\tau L_h)\tilde{e}^n,(-\Delta_N)^{-1}\tilde{e}^{n+1}\right\rangle\\
&+\tau\left\langle(I+\tau L_h)(\phi_{0}+\phi_{1})(\tau L_h)\Delta_N[f_\kappa(U^n)-f_\kappa(\tilde{u}(t_n))],(-\Delta_N)^{-1}\tilde{e}^{n+1}\right\rangle\\
&-\tau\left\langle(I+\tau L_h)\phi_{1}(\tau L_h)\Delta_N[f_\kappa(U^{n-1})-f_\kappa(\tilde{u}(t_{n-1}))],(-\Delta_N)^{-1}\tilde{e}^{n+1}\right\rangle\\
&-\int_0^{\tau}\left\langle(I+\tau L_h)\tilde{e}^{-(\tau-s)L_h}\tilde{R}_{h\tau}^{(2)}(s),(-\Delta_N)^{-1}\tilde{e}^{n+1}\right\rangle\diff s.
\end{align*}
By induction argument, assuming for the numerical error function at the previous time steps $t_{n-1}$ and $t_n$ satisfy
\begin{equation}\label{equation:eq-(3.46)}
\|\tilde{e}^k\|_\infty\leq \frac{1}{2},\quad (k=n,n-1),
\end{equation}
then we have
$
\|U^k\|_\infty=\|\tilde{u}(t_k)+\tilde{e}^k\|_\infty \leq \|\tilde{u}(t_k)\|_\infty +\|\tilde{e}^k\|_\infty\leq \|u_N(t_k)\|_\infty+\frac{1}{2}+\frac{1}{2}\leq B,(k=n,n-1).
$
Based on the above assumptions and $\|\tilde{u}\|_\infty\leq B$, we have
\begin{equation*}
\|f_\kappa(U^k)-f_\kappa(\tilde{u}(t_k))\|_2\leq |3B^2-\kappa-1|\|\tilde{e}^k\|_2:=K\|\tilde{e}^k\|_2,\quad k=n,n-1.
\end{equation*}

We now estimate $R_2$ defined in \eqref{equation:eq-(3.45)}.
By Lemma \ref{lem3} and Cauchy inequality, we obtain
\begin{align}\label{equation:eq-(3.47)}
R_2\leq&\interleave (I+\tau L_h)\phi_{-1}(\tau L_h)\interleave\|\tilde{e}^n\|_{-1,N}\|\tilde{e}^{n+1}\|_{-1,N}\nonumber\\
&+\tau\interleave (I+\tau L_h)(\phi_{0}+\phi_{1})(\tau L_h)\interleave \|f_\kappa(U^n)-f_\kappa(\tilde{u}(t_n))\|_2\|\tilde{e}^{n+1}\|_2\nonumber\\
&+\tau\interleave(I+\tau L_h)\phi_{1}(\tau L_h)\interleave \|f_\kappa(U^{n-1})-f_\kappa(\tilde{u}(t_{n-1}))\|_2\|\tilde{e}^{n+1}\|_2\nonumber\\
&+\int_0^{\tau}\interleave (I+\tau L_h)\tilde{e}^{-(\tau-s)L_h}\interleave\diff s\sup_{t\in(0,\tau)}\|\tilde{R}_{h\tau}^{(2)}(t)\|_{-1,N}\|\tilde{e}^{n+1}\|_{-1,N}\nonumber\\
\leq&\frac{1}{2}\|\tilde{e}^n\|_{-1,N}^2+\frac{1}{2}\|\tilde{e}^{n+1}\|_{-1,N}^2+\frac{3}{2}K\tau(\|\tilde{e}^n\|_2^2+\|\tilde{e}^{n+1}\|_2^2)\nonumber\\
&+\frac{1}{2}K\tau(\|\tilde{e}^{n-1}\|_2^2+\|\tilde{e}^{n+1}\|_2^2)+\frac{\tau}{2}\sup_{t\in(0,\tau)}\|\tilde{R}_{h\tau}^{(2)}(t)\|_{-1,N}^2
+\frac{\tau}{2}\|\tilde{e}^{n+1}\|_{-1,N}^2.
\end{align}
Then, a substitution \eqref{equation:eq-(3.47)} and \eqref{equation:eq-(3.24)} to \eqref{equation:eq-(3.45)}, we obtain
\begin{align}\label{equation:eq-(3.49)}
\frac{1}{2}(\|\tilde{e}^{n+1}\|_{-1,N}^2-\|\tilde{e}^n\|_{-1,N}^2)+(\kappa+1)\tau\|\tilde{e}^{n+1}\|_2^2
&\leq2K\tau\|\tilde{e}^{n+1}\|_2^2+\frac{3}{2}K\tau\|\tilde{e}^n\|_2^2+\frac{1}{2}K\tau\|\tilde{e}^{n-1}\|_2^2\nonumber\\
&+(\frac{C_1}{\gamma_0}+\frac{1}{2})\tau\|\tilde{e}^{n+1}\|_{-1,N}^2+\frac{\tau}{2}\sup_{t\in(0,\tau)}\|\tilde{R}_{h\tau}^{(2)}(t)\|_{-1,N}^2.
\end{align}
Summing the above inequality from 1 to $n$ leads to
\begin{align}\label{equation:eq-(3.51)}
&\|\tilde{e}^{n+1}\|_{-1,N}^2+2(\kappa+1)\tau\sum_{k=1}^{n+1}\|\tilde{e}^k\|_2^2\nonumber\\
&\leq  \|\tilde{e}^1\|_{-1,N}^2+ (\frac{2C_1}{\gamma_0}+1)\tau\sum_{k=1}^{n+1}\|\tilde{e}^k\|_{-1,N}^2+8K\tau\sum_{k=1}^{n+1}\|\tilde{e}^k\|_2^2+\tau\sum_{k=1}^{n}\sup_{t\in(0,\tau)}\|\tilde{R}_{h\tau}^{(2)}(t)\|_{-1,N}^2.
\end{align}
Then using the estimate $\|\tilde{e}^1\|_{-1,N}\leq C(\tau^3+h^m)$ and assuming that $\tau<(\frac{2C_1}{\gamma_0}+1)^{-1}$ and $\gamma_2:=2(\kappa+1)-8K\geq 0$, an application of the discrete Gronwall's inequality gives the  convergence result:
\begin{equation}\label{equation:eq-(3.52)}
\|\tilde{e}^{n+1}\|_{-1,N}+ \left( \gamma_2\tau\sum_{k=1}^{n+1}\|\tilde{e}^k\|_2^2 \right)^{1/2}\leq C_2(\tau^3+h^m),
\end{equation}
where $C_2$ is depends on $C_1, \gamma_0$ and $T$ but independent of $\tau$ and $h$.
Similar to the estimation of \eqref{equation:eq-(3.291)} and \eqref{equation:eq-(3.292)}, we can get
$
\|\tilde{e}^{n+1}\|_\infty \leq \frac{1}{2}$ and $\|U^{n+1}\|_\infty\leq B$
provided that $\tau\leq Ch, m\geq 3$.

Finally, the error estimate \eqref{equation:eq2} can be concluded from \eqref{equation:eq-(3.37)}, \eqref{equation:eq-(3.52)} and the uniform boundedness of $u_{\tau,3}$. This competes the proof of Theorem \ref{thm:th1} (ii).

\section{Numerical experiments}\label{sec:Numer}
In this section, we present numerical experiments to verify the temporal convergence rates in the discrete $H^{-1}$ norm of the proposed ETD1 and ETD2 schemes. We also investigate the phase transition till the steady state and simulate the coarsening dynamics to show the longtime behavior for the NCH equation \eqref{equation:eq-(1.1)}.  Set $\kappa=2$ for the ETD1 scheme and  $\kappa=3$ for the ETD2 scheme in all experiments.
We use the Guass-type function \cite{Du1,zhang2024}
\[
J( x)=\frac{4}{\pi^{d/2}\delta^{d+2}}e^{-\frac{| x|^2}{\delta^2}},\quad x\in \mathbb{R}^d,
\]
where $\delta>0$ is a parameter. Note that $J*1=4/\delta^2$, then $\gamma_0:=\epsilon^2(J\ast 1)-1>0$ is equivalent to $\delta<2 \varepsilon$.
\subsection{Convergence tests}

\begin{exm}[2D test]
We consider the NCH equation \eqref{equation:eq-(1.1)} on $\Omega=(-1,1)^2$ subject to the periodic boundary condition with the initial value $u_0(x,y)=0.05\sin\pi x\sin\pi y$.
\end{exm}
In order to test the time accuracy of the numerical methods, we calculate the $H^{-1}$-norm error using
\[
e(\tau):=\left[ h^2\sum_{i,j=0}^{N-1}\left((-\Delta_N)^{-\frac{1}{2}}\left(U_{i,j}^{N_t}(\tau,h)-U_{i,j}^{2N_t}(\tau/2,h)\right)\right) \right]^\frac{1}{2}
\]
and the time convergence order by $\textrm{Rate}=\ln[e(\tau)/e(\tau/2)]/\ln 2$.

In this simulation $N=256$ and  $T=0.5$. We compute the numerical solutions for ETD1 and ETD2 schemes with various time step sizes $\tau=0.005\times 2^{-k}, k=0,1,\cdots,8$. The reference solutions is taken as the approximated solution obtained with a quite small time step size $\tau=0.001\times 2^{-8}$. The errors and convergence rates in the discrete $H^{-1}$ norm are showed in Table \ref{table 1} and Table \ref{table 11}, which shows the first order of accuracy for the ETD1 scheme and the second order of accuracy for the ETD2 scheme with different $\delta$ and $\varepsilon$, as expected.

\begin{table}[!htb]
	\caption{\label{table 1}  The $H^{-1}$ errors and convergence rates of the ETD1 scheme at time $T=0.5$.}
	\begin{center}
		\begin{tabular}{|l c| cc |cc|}
        \hline
      \multicolumn{2}{|c|}{ETD1}&\multicolumn{2}{c|}{$\delta^2=\varepsilon^2=0.1$} &\multicolumn{2}{c|}{$\delta^2=\varepsilon^2=0.01$}\\
      \hline
			 &$\tau=0.005$ & {$H^{-1}$ error } & {Rate} & {$H^{-1}$ error } & {Rate} \\ \hline
			 & $\tau$ &   5.1108E-05  &  --&7.2409E-05 & --\\
		     & $\tau/2$ & 2.1838E-05  & 1.2267 &2.1994E-05&1.7191\\
			 & $\tau/4$  & 1.0061E-05  & 1.1181& 8.7574E-06&1.3285  \\
			 & $\tau/8$  & 4.8154E-06  &  1.0630&3.9203E-06&1.1595 \\
			 & $\tau/16$ & 2.3443E-06  & 1.0385& 1.8485E-06  & 1.0846  \\
             & $\tau/32$ & 1.1456E-06  & 1.0331& 8.8928E-07  & 1.0556  \\
             & $\tau/64$ & 5.5526E-07 & 1.0448& 4.2773E-07 & 1.0560  \\
             & $\tau/128$ & 2.6235E-07  & 1.0817& 2.0132E-07  & 1.0872  \\
             & $\tau/256$ & 1.1645E-07  & 1.1718& 8.9191E-08  & 1.1745  \\ \hline
		\end{tabular}
	\end{center}
\end{table}

\begin{table}[!htb]
	\caption{\label{table 11}  The $H^{-1}$ errors and convergence rates of the ETD2 scheme at time $T=0.5$.}
	\begin{center}
		\begin{tabular}{|l c| cc |cc|}
        \hline
      \multicolumn{2}{|c|}{ETD2}&\multicolumn{2}{c|}{$\delta^2=\varepsilon^2=0.1$} &\multicolumn{2}{c|}{$\delta^2=\varepsilon^2=0.01$}\\
      \hline
			 &$\tau=0.005$ & {$H^{-1}$ error } & {Rate} & {$H^{-1}$ error } & {Rate} \\ \hline
			 & $\tau$ &   1.1417E-05  &  --&1.8032E-06 & --\\
		     & $\tau/2$ & 2.5795E-06  & 2.1461 &3.9009E-07&2.2087\\
			 & $\tau/4$  & 6.2025E-07  & 2.0562& 8.8718E-08&2.1365  \\
			 & $\tau/8$  & 1.5266E-07  &  2.0225&2.0981E-08&2.0801 \\
			 & $\tau/16$ & 3.7907E-08  & 2.0098& 5.0878E-09&2.0040  \\
             & $\tau/32$ & 9.4434E-09  & 2.0051& 1.2512E-09&2.0238  \\
             & $\tau/64$ & 2.3530E-09 & 2.0048&3.0946E-10&2.0155  \\
             & $\tau/128$ & 5.8343E-10  & 2.0119&7.6282E-11& 2.0203  \\
             & $\tau/256$ & 1.4140E-10  & 2.0448&1.8444E-11&2.0482  \\ \hline
		\end{tabular}
	\end{center}
\end{table}

\begin{exm}[3D test]
In this test, we calculate the errors and convergence rates by ETD2 scheme for the NCH equation on the computational domain $(-1,1)^3\times(0,T]$. We choose the same $N$ and $\tau$ as 2D.
\end{exm}

We consider the smooth initial condition $u_0(x,y,z)=0.05\sin\pi x\sin\pi y\sin\pi z$, the computational results are presented in Table \ref{table 111}. From Table \ref{table 111} one can see that time accuracy is second order, which is onsistent with our theoretical predictions.

\begin{table}[!htb]
	\caption{\label{table 111}  The $H^{-1}$ errors and convergence rates of the ETD2 scheme at time $T=0.5$.}
	\begin{center}
		\begin{tabular}{|l c| cc |cc|}
        \hline
      \multicolumn{2}{|c|}{ETD2}&\multicolumn{2}{c|}{$\varepsilon=0.2, \delta=0.1$} &\multicolumn{2}{c|}{$\varepsilon=0.2, \delta=0.2$}\\
      \hline
			 &$\tau=0.005$ & {$H^{-1}$ error } & {Rate} & {$H^{-1}$ error } & {Rate} \\ \hline
			 & $\tau$ &   1.8586E-05  &  --&2.8680E-05 & --\\
		     & $\tau/2$ & 3.5861E-06  & 2.3737 &2.4139E-06&3.5706\\
			 & $\tau/4$  & 8.1176E-07  & 2.1433& 2.6908E-07&3.1652  \\
			 & $\tau/8$  & 1.9593E-07  &  2.0507&4.4968E-08&2.5811 \\
			 & $\tau/16$ & 4.8344E-08  & 2.0189& 9.7783E-09&2.2012  \\
             & $\tau/32$ & 1.2017E-08  & 2.0082& 2.3470E-09&2.0587  \\
             & $\tau/64$ & 2.9918E-09 & 2.0060&5.7920E-10&2.0187  \\
             & $\tau/128$ & 7.4157E-10  & 2.0123&1.4335E-10& 2.0146  \\
             & $\tau/256$ & 1.7970E-10  & 2.0450&3.4806E-11&2.0421  \\ \hline
		\end{tabular}
	\end{center}
\end{table}

\subsection{Interfaces in the steady states}
\begin{figure*}[!htb]
	\centering
    {\includegraphics[width=0.3\textwidth]{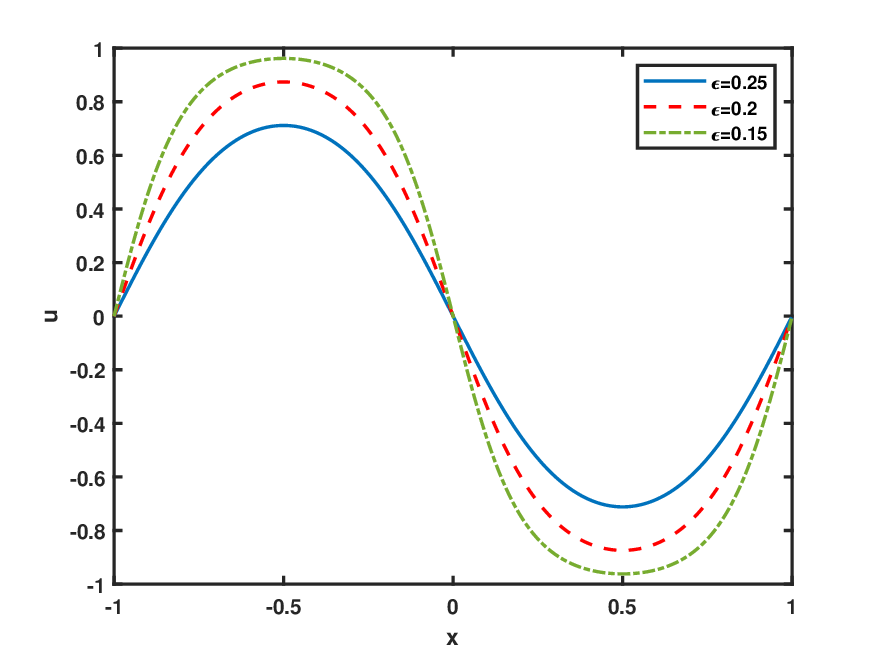}}
    {\includegraphics[width=0.3\textwidth]{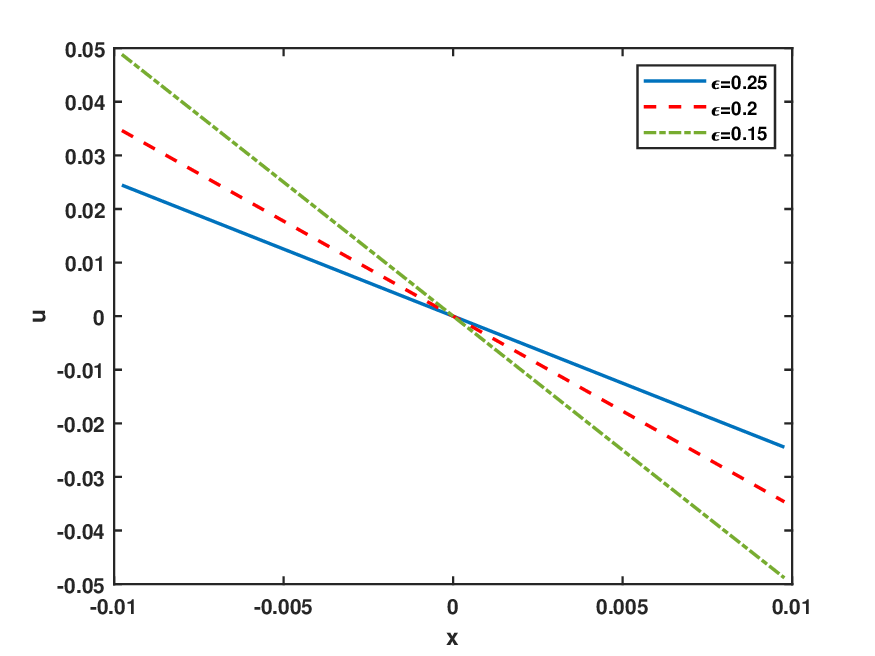}}
    {\includegraphics[width=0.3\textwidth]{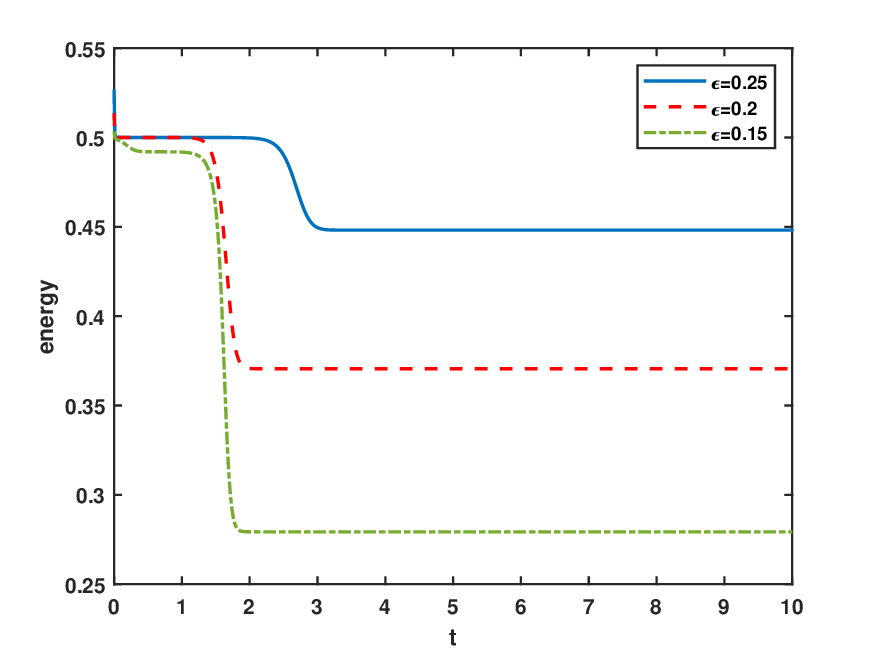}}
    \subfigure[Steady states]{\includegraphics[width=0.3\textwidth]{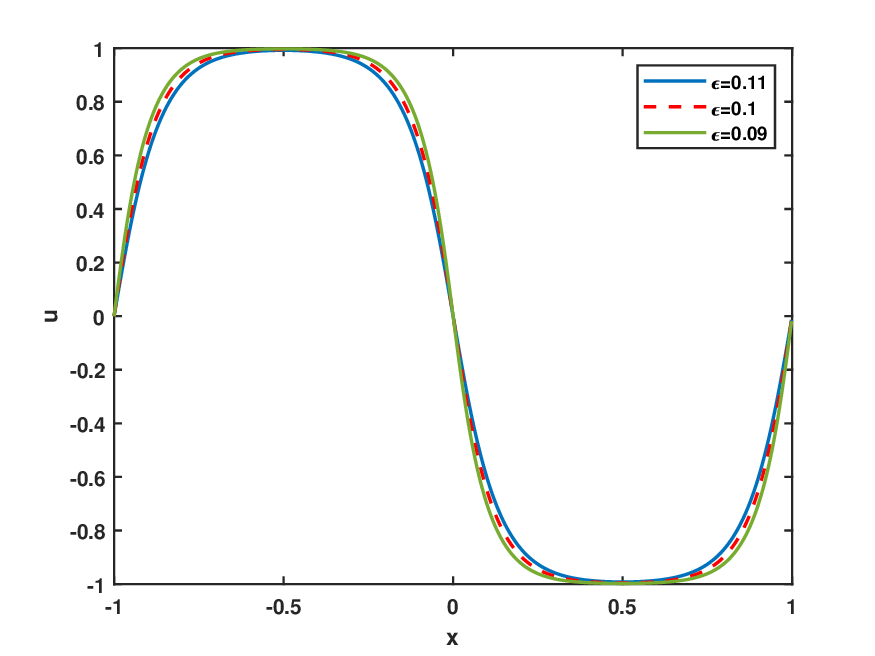}}
    \subfigure[Phases around one interface]{\includegraphics[width=0.3\textwidth]{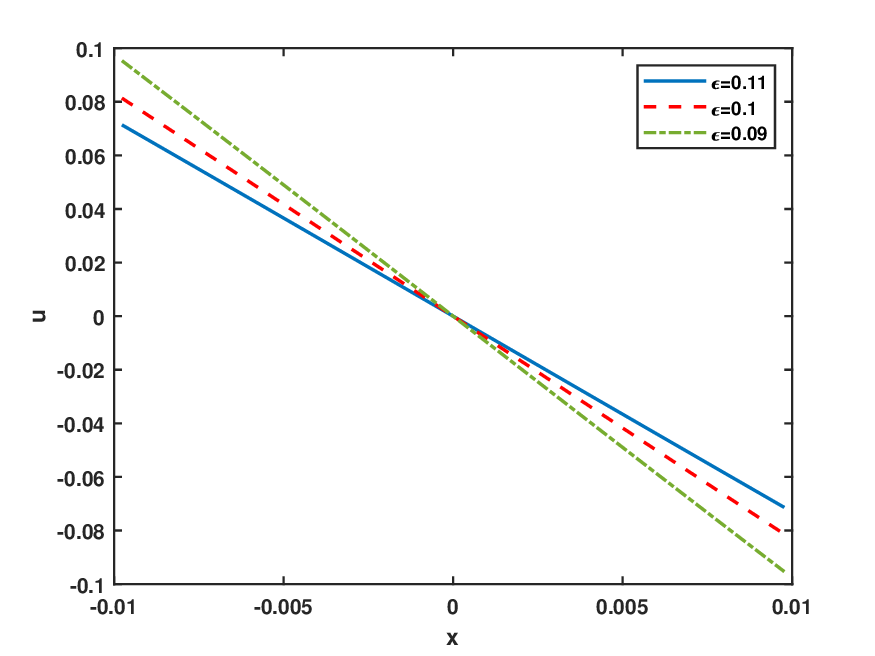}}
    \subfigure[Evolutions of energy]{\includegraphics[width=0.3\textwidth]{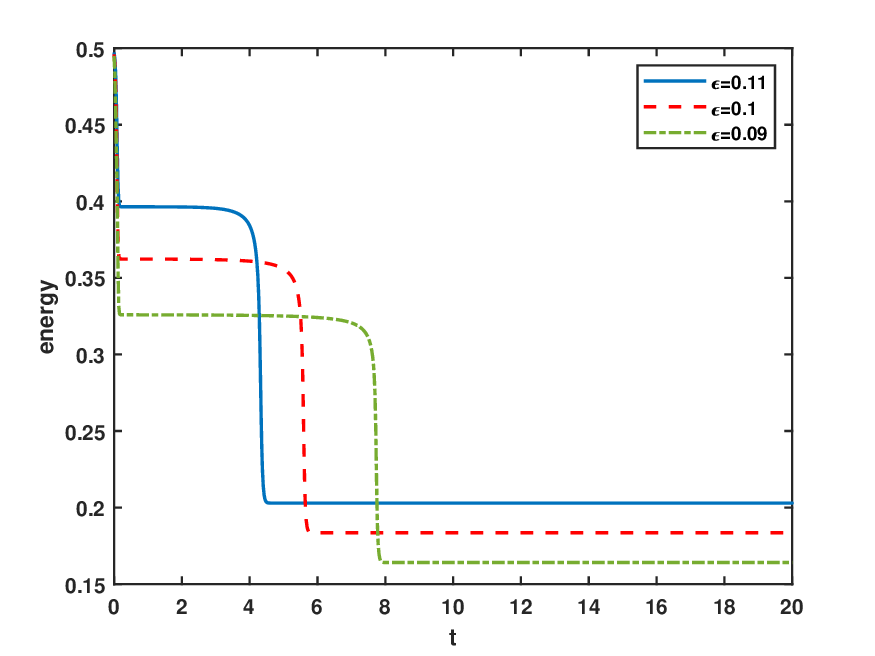}}
	\caption{Numerical simulation with $\delta=0.1$ for different values of $\varepsilon$.}
    \label{Figure1}
    \vspace{1cm}
    \subfigure[Steady states]{\includegraphics[width=0.3\textwidth]{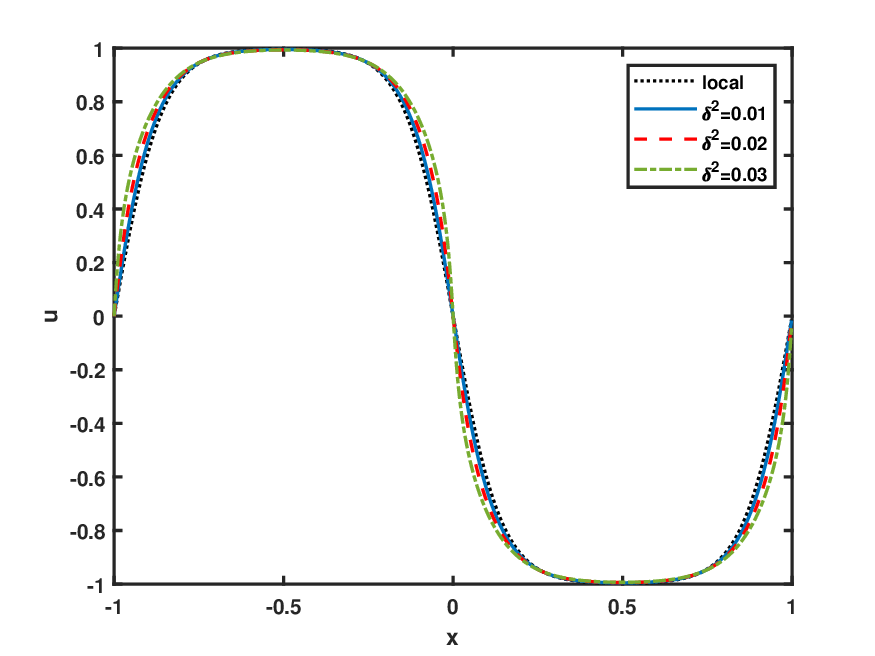}}
    \subfigure[Phases around one interface]{\includegraphics[width=0.3\textwidth]{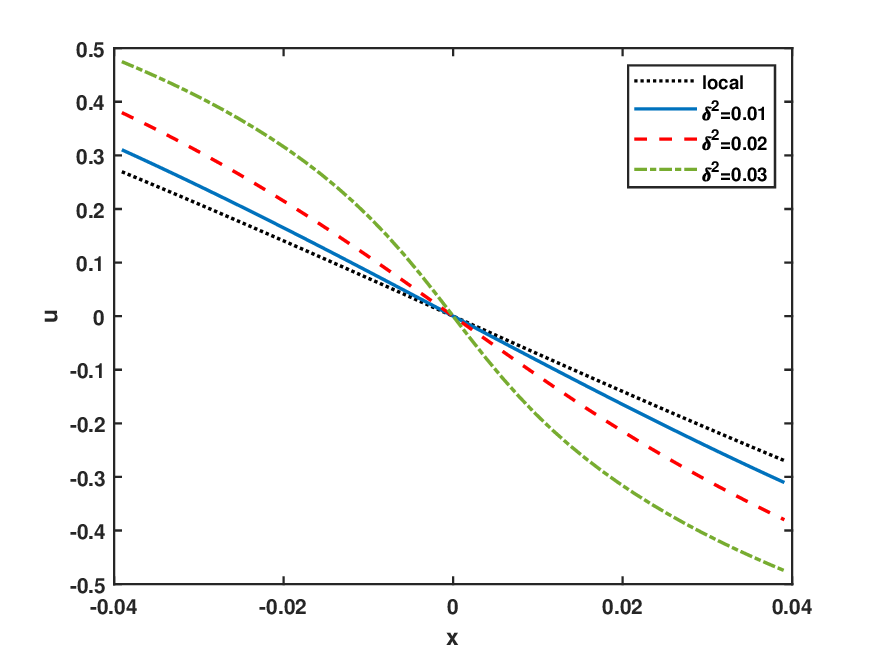}}
    \subfigure[Evolutions of energy]{\includegraphics[width=0.3\textwidth]{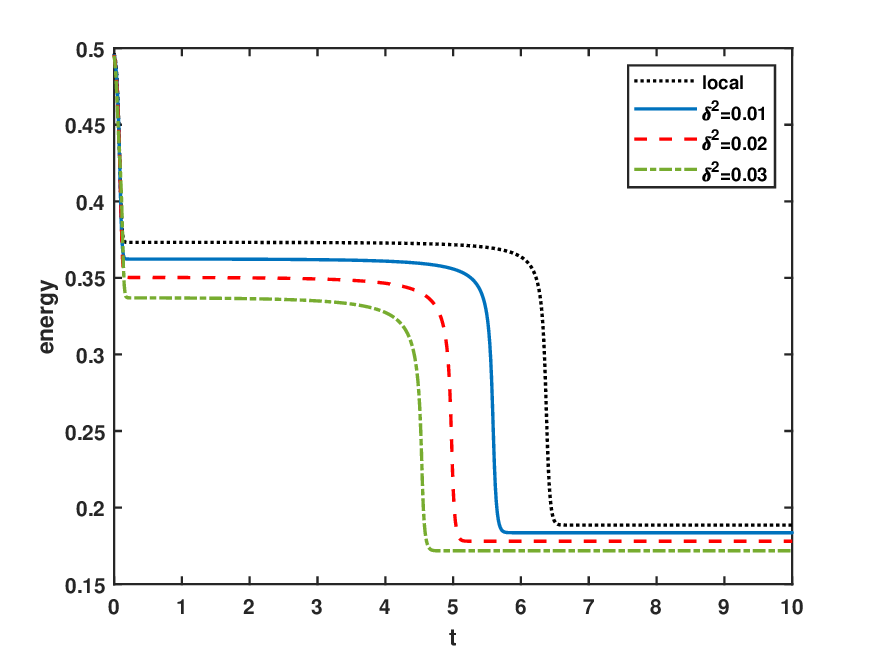}}
	\caption{Numerical simulation with $\varepsilon=0.1$ for different values of $\delta$.}
    \label{Figure2}
\end{figure*}

We now simulate the shapes of the interfaces formed in the steady states by the NCH equation \eqref{equation:eq-(1.1)} under the ETD2 scheme in one-dimensional case with various $\varepsilon$ and $\delta$.
\begin{exm}[1D problem]
Let $\Omega=(-1,1)$, $u_0(x)=0.1(\sin 2\pi x+\sin 3\pi x)$, $N=1024$ and $\tau=0.0001$.
\end{exm}
In Figure \ref{Figure1}, we fix $\delta=0.1$ and choose $\varepsilon$ as 0.25, 0.2, 0.15, 0.11, 0.1 and 0.09. We plot the numerical solutions at the steady states and the curve of discrete energy evolution. The Figure shows that the numerical solutions reach the steady state at $T=10$ and $T=20$, respectively. The energy stability of the ETD2 can be seen in the last column. It is observed that the time required to reach the steady state increases as $\varepsilon$ decreases and the interface width becomes sharper for smaller $\varepsilon$.

In Figure \ref{Figure2}, we choose $\varepsilon=0.1$, and $\delta^2$ as 0.03, 0.02, 0.01 respectively. Furthermore, we simulate the local CH equation for comparison. From Figure \ref{Figure2}, we observe the energy curve remain unchanged after $T=10$ and the numerical solutions reach the steady state in this time. As $\delta$ decreases, the interface turns flat, and the phase transition process is close to the local one for all cases.

\subsection{Coarsening dynamics and energy evolution}
\begin{exm}[2D problem]
We now simulate the long time behavior of the NCH equation \eqref{equation:eq-(1.1)} by using the ETD2 scheme in $\Omega=(-2\pi,2\pi)\times (-2\pi,2\pi)$ with a random initial data ranging from -0.1 to 0.1. We set time step size $\tau=0.01$ and choose $N=512$.
\end{exm}
\begin{figure*}[!htb]
	\centering
    \subfigure[T=1]{\includegraphics[width=0.33\textwidth]{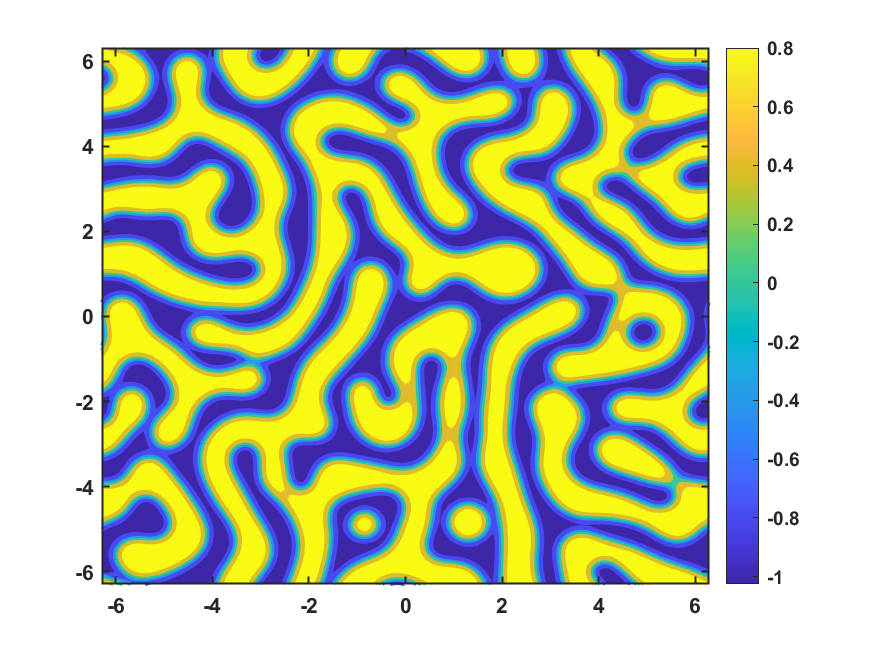}}
    \subfigure[T=10]{\includegraphics[width=0.33\textwidth]{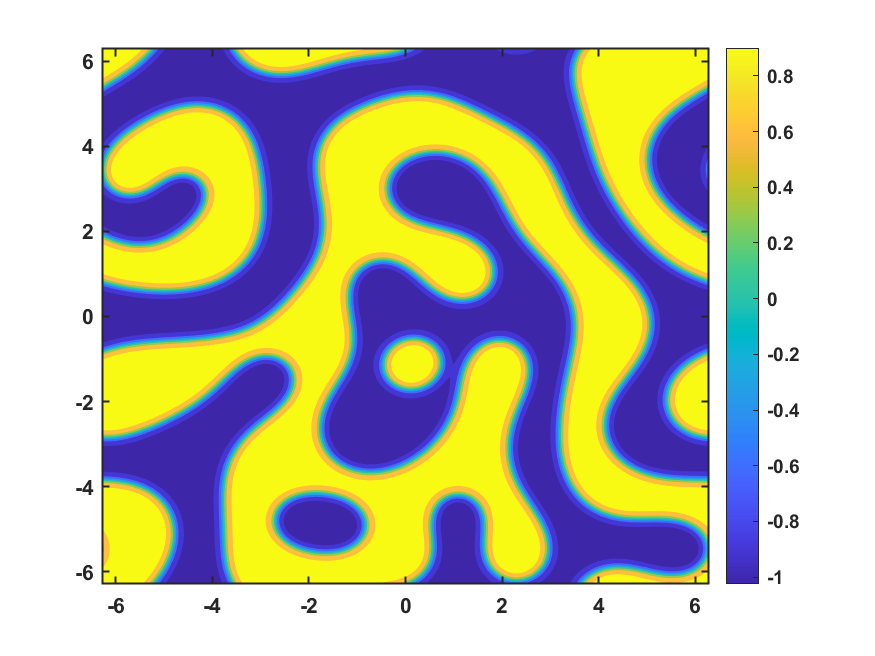}}
    \subfigure[T=100]{\includegraphics[width=0.33\textwidth]{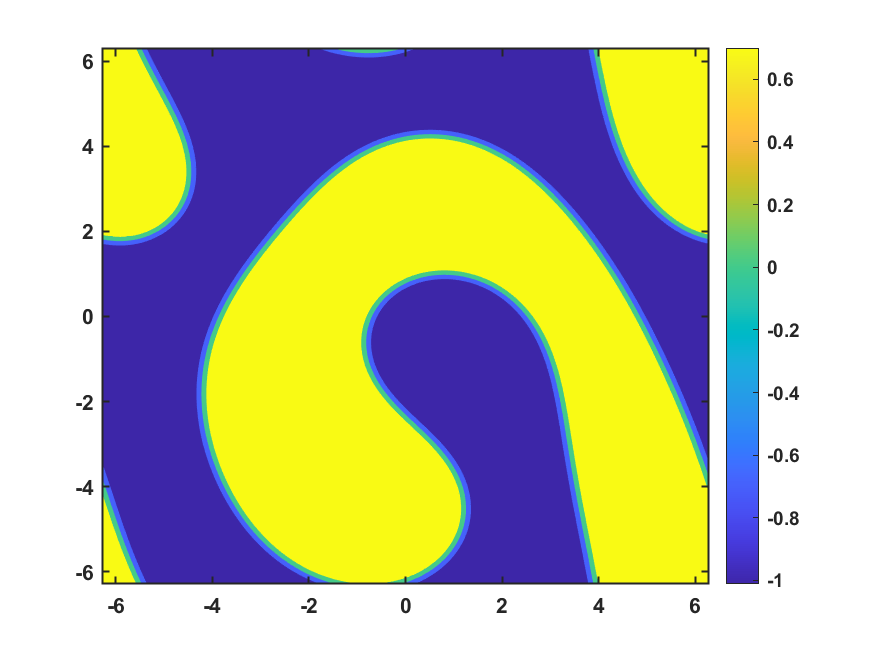}}
    \subfigure[T=400]{\includegraphics[width=0.33\textwidth]{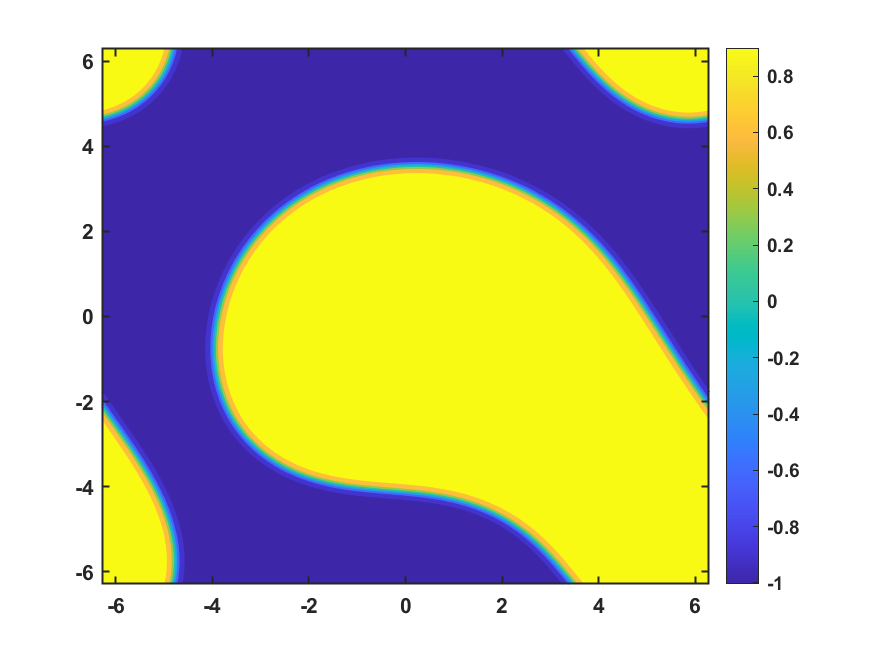}}
    \subfigure[T=1200]{\includegraphics[width=0.33\textwidth]{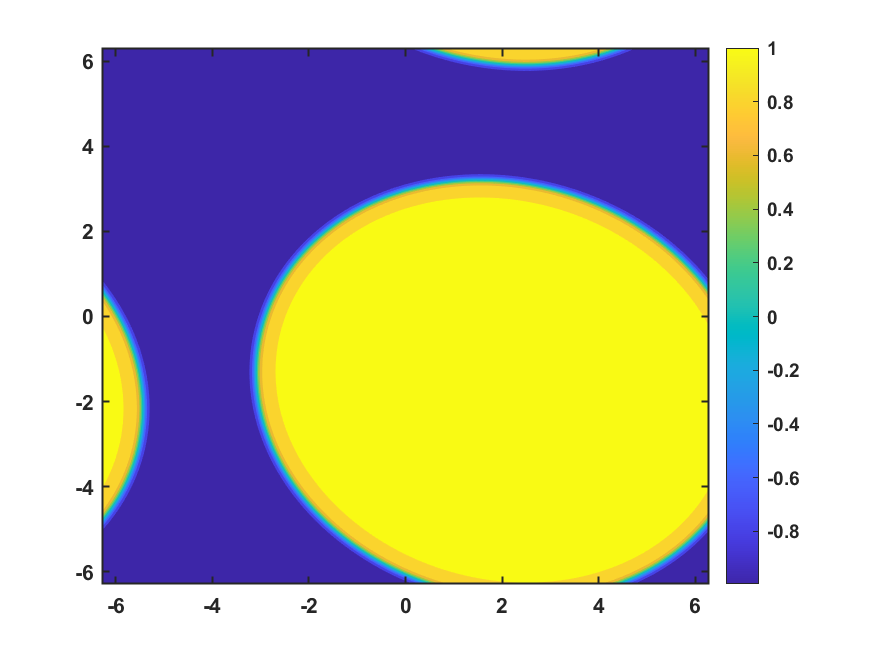}}
    \subfigure[T=2000]{\includegraphics[width=0.33\textwidth]{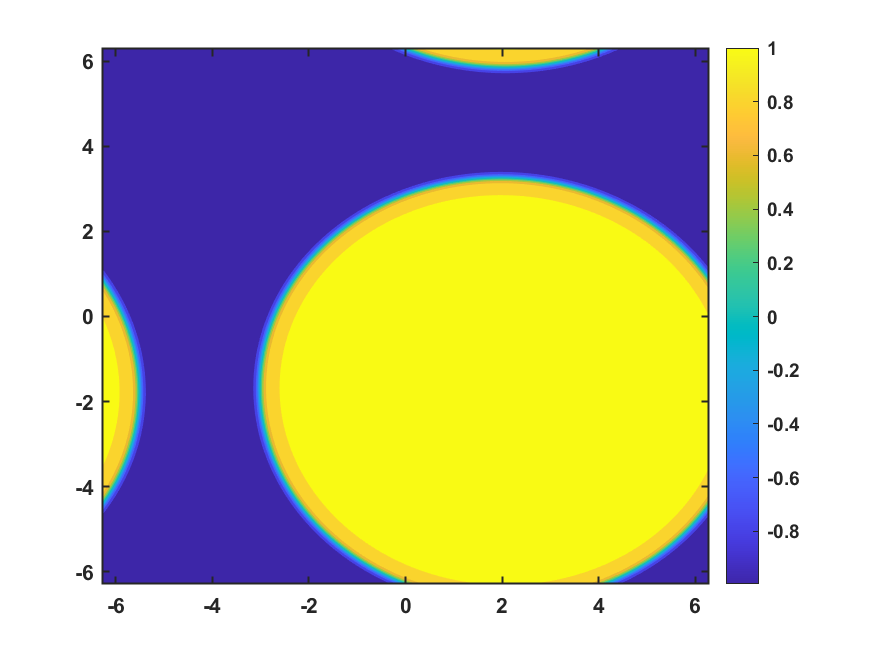}}
	\caption{Numerical results at $T=1,10,100,400,1200,2000$.}
    \label{Figure3}
    \vspace{1cm}
	\centering
     \subfigure[$\delta=0.05, N=512$]{\includegraphics[width=0.45\textwidth]{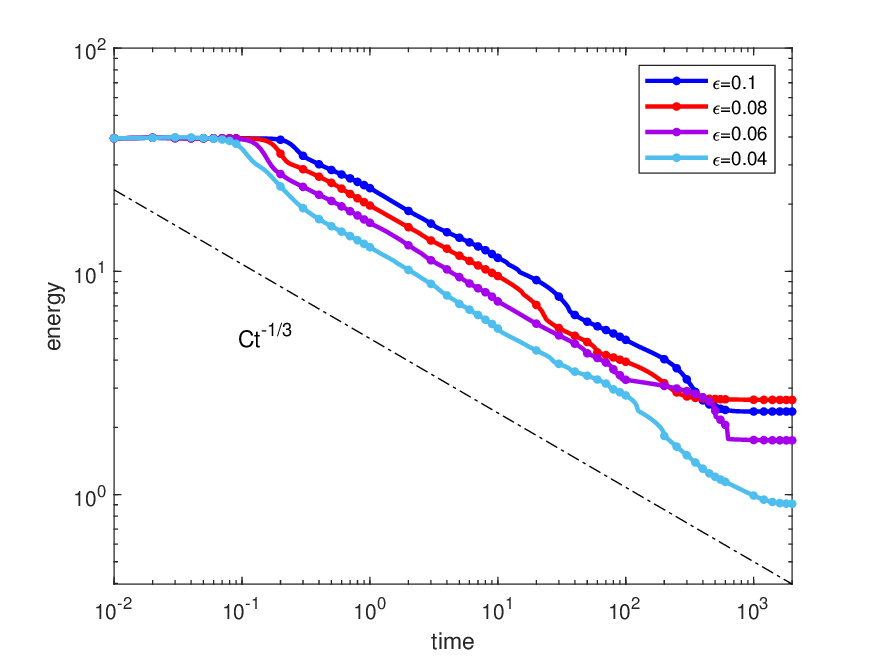}\label{fig: figure1}}
     \subfigure[$\varepsilon=0.1, N=512$]{\includegraphics[width=0.45\textwidth]{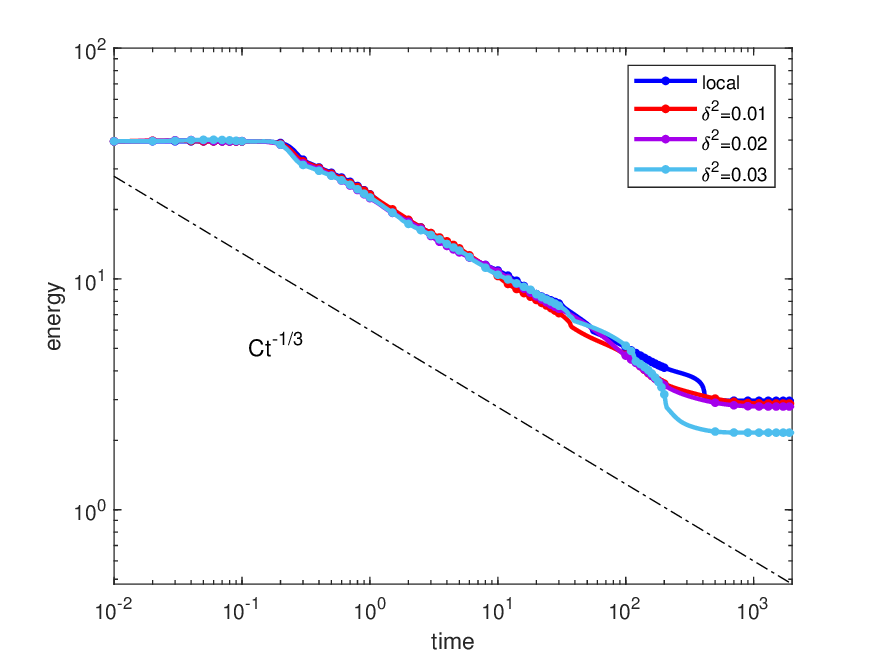}\label{fig: figure2}}
	\caption{The temporal evolution of the energy.}
    \label{fig: figure}
\end{figure*}

Let $\varepsilon=\delta=0.09$. Figure \ref{Figure3} shows the coarsening dynamics of numerical solutions, from which one can observe that the dynamic evolves from the initial disorder state to the ordered states and reach the steady state around $T=2000$.
\begin{figure*}[!htb]
	\centering
    \subfigure[T=0]{\includegraphics[width=0.33\textwidth]{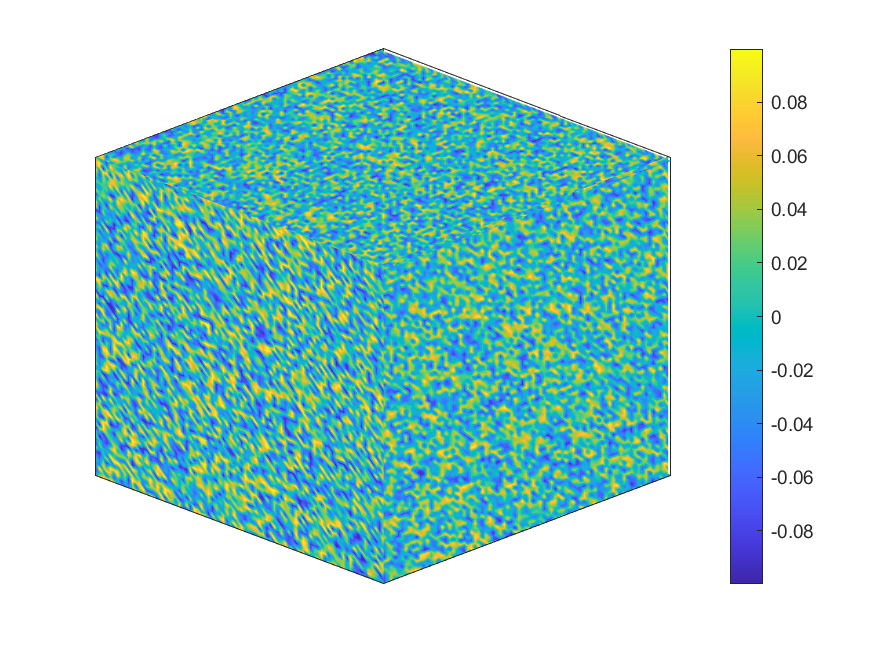}}
    \subfigure[T=50]{\includegraphics[width=0.33\textwidth]{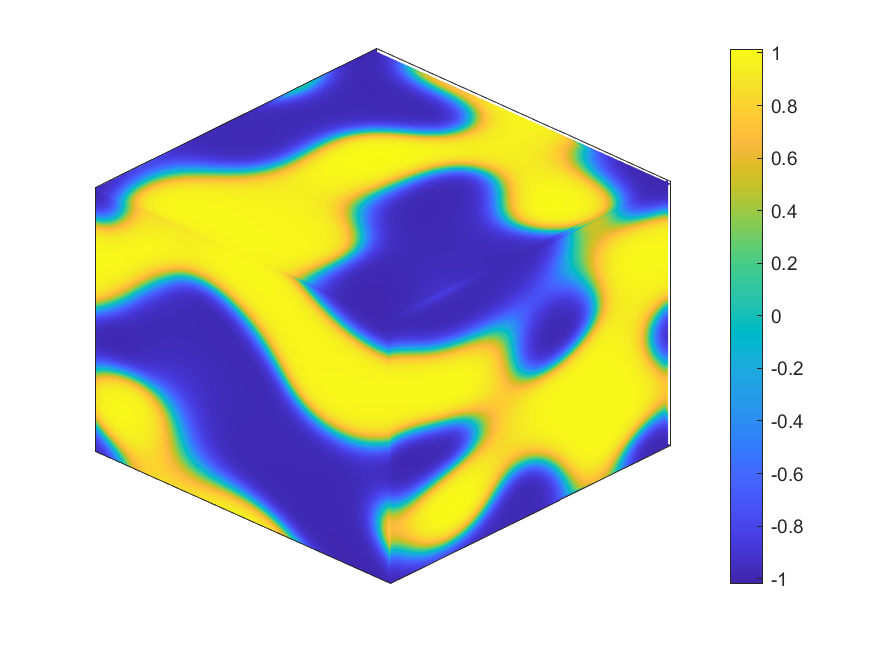}}
    \subfigure[T=200]{\includegraphics[width=0.33\textwidth]{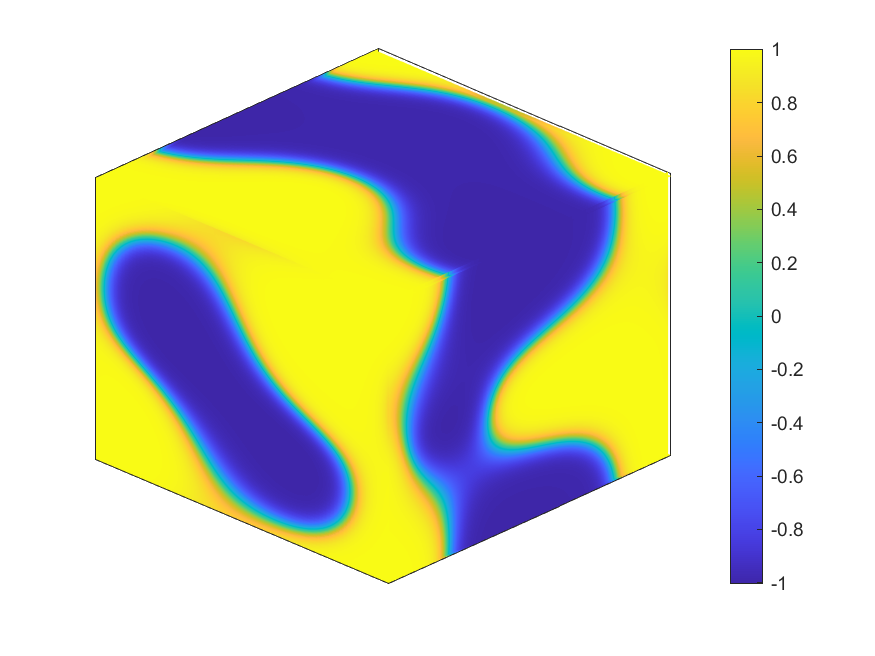}}
    \subfigure[T=500]{\includegraphics[width=0.33\textwidth]{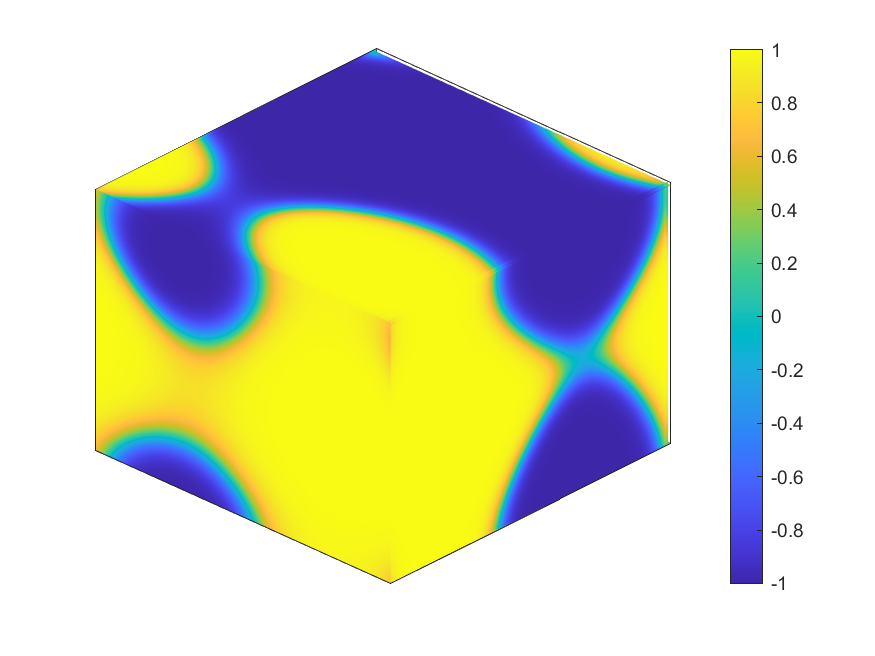}}
    \subfigure[T=700]{\includegraphics[width=0.33\textwidth]{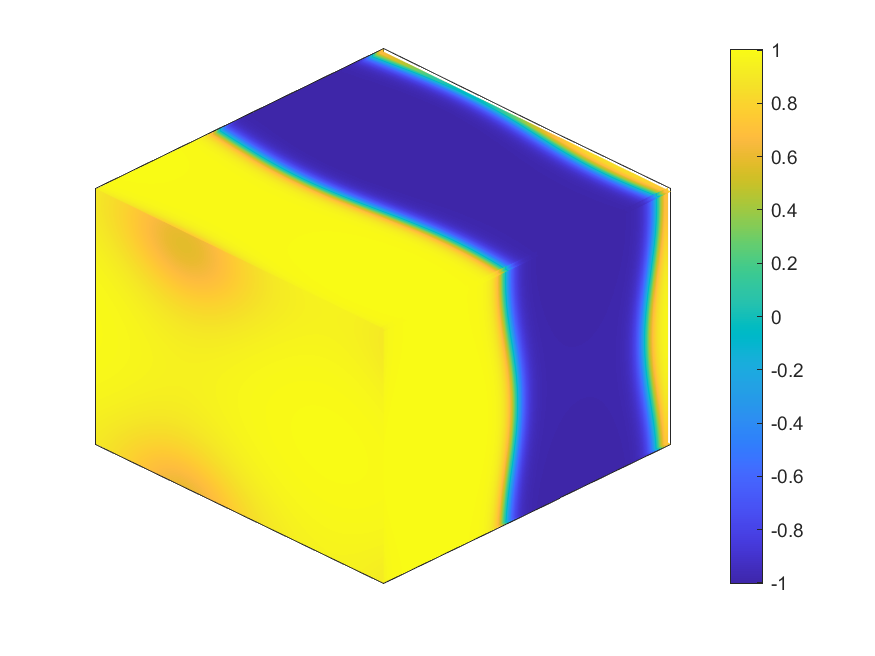}}
    \subfigure[T=1000]{\includegraphics[width=0.33\textwidth]{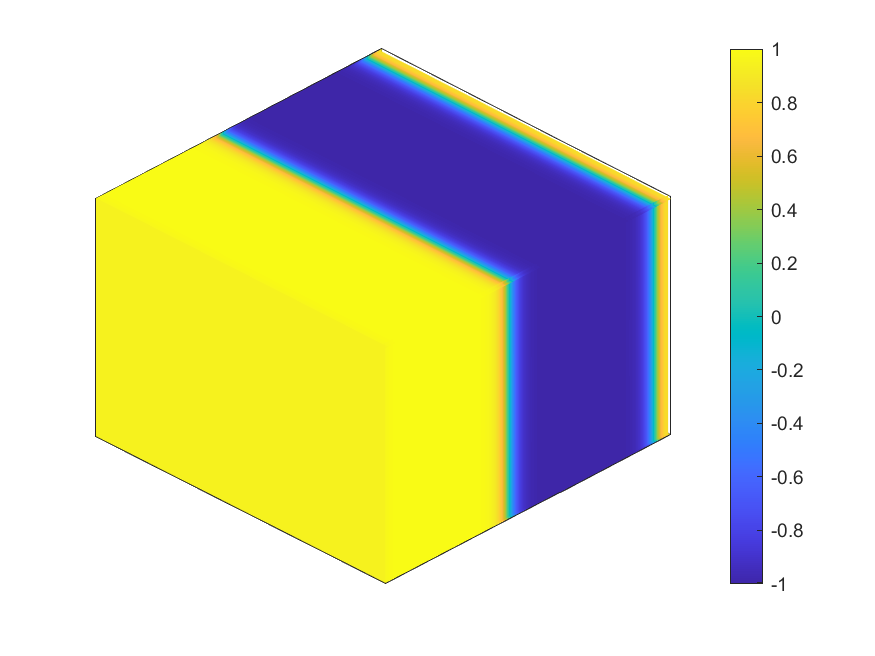}}
	\caption{Numerical results at $T=0,50,200,500,700,1000$.}
    \label{Figure4}
    \vspace{1cm}
	\centering
     \subfigure{\includegraphics[width=0.5\textwidth]{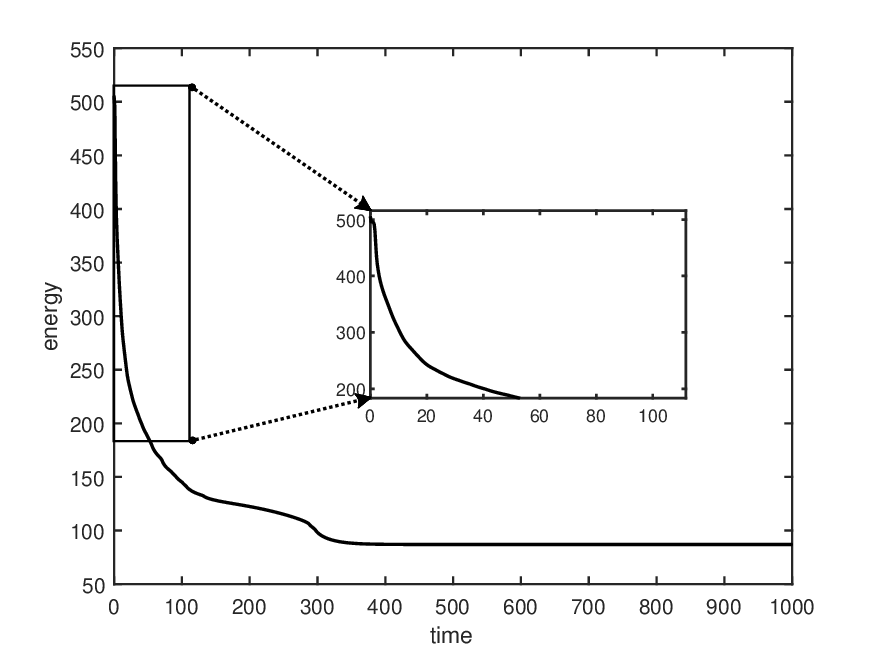}}
	\caption{The temporal evolution of the energy.}
    \label{Figure5}
\end{figure*}

Let $\delta=0.05$ and $\varepsilon=0.1,0.08,0.06,0.04$, respectively. Figure \ref{fig: figure1} plots the energy evolution curves, and we can observe the energy decay rates comply with the $t^{-\frac{1}{3}}$ power law for all cases.
This is consistent to the local CH equation \cite{dai,ju2015}. We also plot the curves of the energy for $\varepsilon=0.1$ with different $\delta$ in Figure \ref{fig: figure2}. Table \ref{table 2} presents the linear fitting coefficients $m_e$ and $b_e$ for $\varepsilon$ decreasing from $0.1$ to $0.04$ with the fitting function $E(t)=b_et^{m_e}$.

\begin{table}[!htb]
	\caption{\label{table 2}  The linear fitting coefficients for the case $\delta=0.05$.}
	\begin{center}
		\begin{tabular}{l c c c c cc cc}
        \hline
			 &$\varepsilon$ &  0.1 & 0.09 & 0.08 & 0.07& 0.06&0.05 & 0.04 \\ \hline
			 & $m_e$ &    -0.314  &-0.314 & -0.325& -0.331&-0.303&-0.320 & -0.341\\
		     & $b_e$ &  21.08  &19.49 & 18.15&17.13 &15.04&13.46 &12.64\\ \hline
		\end{tabular}
	\end{center}
\end{table}

\begin{exm}[3D problem]
We perform the coarsening dynamics of the numerical solutions in 3D on the domain $\Omega=(-2\pi,2\pi)^3$ with initial data $u_0(x,y,z)=-0.1+0.2\cdot \text{rand}(x,y,z)$. Set $N=80$, $\tau=0.01$ and $\varepsilon=\delta=0.3$.
\end{exm}

In Figure \ref{Figure4}, we plot the contours for numerical solutions at different times. It is shown the numerical solution reaches equilibrium state at around $T=1000$ and the energy decay curve is shown in Figure \ref{Figure5}.

\section{Conclusion remarks}\label{sec:Conclu}
Because of the non-locality and nonlinearity, it is very challenging to prove the convergence of numerical methods for NCH equation.
In this paper, we provide a detailed convergence analysis of the ETD1 and ETD2 schemes for the NCH equation, where the Fourier spectral collocation method is used for the spatial discretizations. Due to the lack of the higher-order diffusion term, we adopt $(-\Delta_N)^{-1}\tilde{e}^{n+1}$ as the test function rather than $\tilde{e}^{n+1}$ in the classic CH equation for error equations. The optimal convergence rates in discrete $H^{-1}$ norm have been presented by performing the high order consistency analysis. In addition, we also obtain $\ell^\infty$ bound of the numerical solutions under some moderate constraints on the space-time step sizes. Finally, we verify the convergence in time of the proposed schemes numerically and simulate the coarsening dynamics to show the long time behavior and present the power law for the energy decay for the NCH equation.

\section*{Appendix A: Fourier spectral collocation Approximations}\label{sec:Appendix A}
We introduce some notations and useful properties of Fourier spectral collocation approximations for space local linear operators and nonlocal operator in two dimension. Define the index sets
\begin{align*}
    S_h=\{(i,j)\in \mathbb{Z}^2|1\leq i,j\leq N\},\,
    \widehat{S}_h=\left\{ (k,l)\in \mathbb{Z}^2|-\frac{N}{2}+1\leq i,j\leq \frac{N}{2} \right\}.
\end{align*}
The space of all the $\Omega_h$-periodic grid functions is denoted by $\mathcal{M}_h$.
For any $f, g\in \mathcal{M}_h$, 
the discrete $\ell^2$ inner product $\langle\cdot,\cdot\rangle$ and norm $\|\cdot\|_2$ and discrete $\ell^\infty$ norm $\|\cdot\|_\infty$ are, respectively, defined by
\begin{equation*}
\langle f,g\rangle=h^2\sum_{(i,j)\in S_h}f_{ij}g_{ij},\quad \|f\|_2=\sqrt{\langle f,f\rangle},\quad
 \|f\|_\infty=\max_{(i,j)\in S_h}|f_{ij}|.
\end{equation*}

For a function $f\in \mathcal{M}_h$, we have discrete Fourier transform \cite{wang}
\begin{equation}
f_{ij}=\sum_{(k,l)\in \hat{S}_h}\hat{f}_{kl}\exp(\textrm{i}\frac{k\pi}{X}x_i)\exp(\textrm{i}\frac{l\pi}{X}y_j),\quad (i,j)\in S_h.
\end{equation}
where
\begin{equation}\label{equation:eq-2}
\hat{f}_{kl}=\frac{1}{N^2}\sum_{(i,j)\in S_h}f_{ij}\exp(-\textrm{i}\frac{k\pi}{X}x_i)\exp(-\textrm{i}\frac{l\pi}{X}y_j),\quad (k,l)\in \widehat{S}_h,
\end{equation}
are the discrete Fourier coefficients.
The first and second order derivatives of $f$ in the $x$ direction can be defined as
\[
D_xf_{ij}=\sum_{(k,l)\in \hat{S}_h}(\textrm{i}\frac{k\pi}{X})\hat{f}_{kl}\exp(\textrm{i}\frac{k\pi}{X}x_i)\exp(\textrm{i}\frac{l\pi}{X}y_j),\quad
D_x^2f_{ij}=\sum_{(k,l)\in \hat{S}_h}(\textrm{i}\frac{k\pi}{X})^2\hat{f}_{kl}\exp(\textrm{i}\frac{k\pi}{X}x_i)\exp(\textrm{i}\frac{l\pi}{X}y_j).
\]
The differentiation operators in the $y$ direction can be defined in the same way.
For any $f\in\mathcal{M}_h$, the discrete gradient and Laplace operators are given by
\begin{equation}\label{equation:eq-3}
\nabla_Nf=\begin{pmatrix}D_xf\\D_yf\end{pmatrix},\quad \Delta_Nf=D_x^2f+D_y^2f.
\end{equation}

For any $f,g\in\mathcal{M}_h, \textbf{\textit{g}}\in \mathcal{M}_h\times\mathcal{M}_h$, we have the following summation-by-parts formulas \cite{Ju}
\[
\langle f,\nabla_N\cdot \textbf{\textit{g}}\rangle=-\langle\nabla_Nf,\textbf{\textit{g}}\rangle,\quad \langle f,\Delta_Ng\rangle= -\langle\nabla_Nf,\nabla_Ng\rangle=\langle\Delta_Nf,g\rangle.
\]
For any $f\in\mathcal{M}_h$, we call $\bar{f}:=\frac{1}{4X^2}\langle f, 1\rangle$ the mean value of $f$. By noticing \eqref{equation:eq-mass} and assuming that the mean of $u_0$ is zero, we only need to consider the zero-mean grid functions, i.e.,
\[
\mathcal{M}_h^0=\{v\in\mathcal{M}_h|\langle v,1\rangle=0\}=\{v\in\mathcal{M}_h|\hat{v}_{00}=0\},
\]
where $\hat{v}_{00}=\frac{1}{N^2}\sum_{(i, j)\in S_h}v_{ij}=\frac{1}{4X^2}\langle v,1\rangle$.
Similar to the continuous case, $(-\Delta_N)$ is self-adjoint and positive definite on $\mathcal{M}_h^0$ and thus $(-\Delta_N)^{-1}$ is well-defined and is also self-adjoint and positive definite.
Then for any $f, g\in \mathcal{M}_h^0$, we can define the discrete $H_h^{-1}$ inner produce and the discrete $H_h^{-1}$ norm as
\begin{equation*}
\langle f,g\rangle_{-1,N}:=\langle f,(-\Delta_N)^{-1}g\rangle=\langle (-\Delta_N)^{-\frac{1}{2}}f,(-\Delta_N)^{-\frac{1}{2}}g\rangle, \quad
\|f\|_{-1,N}:=\sqrt{\langle f,f \rangle_{-1,N}}=\|(-\Delta_N)^{-\frac{1}{2}}f\|_2.
\end{equation*}

For any $f,\phi\in \mathcal{M}_h$, the discrete convolution $f\circledast \phi\in \mathcal{M}_h$ can be defined componentwise \cite{Du1} by
\[
(f\circledast \phi)_{ij}=h^2\sum_{(m,n)\in S_h}f_{i-m,j-n}\phi_{mn},\quad (i,j)\in S_h.
\]
Given a kernel function $J$ satisfying the conditions,  then for any $f\in\mathcal{M}_h$, the discrete form of nonlocal diffusion operator $\mathcal{L}$  can be defined by $\mathcal{L}_N=\mathcal{F}^{-1}\hat{\mathcal{L}}_N\mathcal{F}$ with
\begin{equation}\label{equation:eq-4}
\hat{\mathcal{L}}_N\hat{f}_{kl}=\lambda_{kl}\hat{f}_{kl},\quad (k,l)\in \hat{S}_h
\end{equation}
where $\mathcal{F}$ is a discrete Fourier transform matrix and
\[
\lambda_{kl}=h^2\sum_{(i,j)\in S_h}J(x_i,y_j)\left(1-\exp(-\textrm{i}\frac{k\pi}{X}x_i)\exp(-\textrm{i}\frac{l\pi}{X}y_j)\right).
\]

\bibliographystyle{alpha}
\bibliography{NCHETD}
\end{document}